\documentclass[10pt,a4paper,twoside]{article}
\usepackage{amsfonts}
\usepackage{amssymb}
\usepackage{amsthm}
\usepackage{latexsym}
\usepackage{amsmath}
\usepackage{amssymb}

\newtheorem{thm}{Theorem}[section]

\newtheorem{lem}[thm]{Lemma}
\newtheorem{remark}[thm]{Remark}
\newtheorem{cor}[thm]{Corollary}
\newtheorem{prop}[thm]{Proposition}

\newtheorem{definition}[thm]{Def\mbox{}inition}

\newcommand{\R}{\mathbb{R}}
\newcommand{\Rr}{\mathcal{R}}
\newcommand{\C}{\mathbb{C}}

\newcommand{\Q}{\mathbb{Q}}
\newcommand{\Z}{\mathbb{Z}}
\newcommand{\Y}{\mathcal{O}}

\newcommand{\s}{\mathfrak{l}}
\newcommand{\f}{\mathfrak{f}}
\newcommand{\p}{\mathfrak{p}}

\newcommand{\Mm}{\mathcal{M}}
\newcommand{\Dd}{\mathcal{D}}
\newcommand{\Cc}{\mathcal{C}}
\newcommand{\Uu}{\mathcal{U}}
\newcommand{\Pol}{\mathcal{P}ol}

\newcommand{\Ee}{\mathcal{E}}

\newcommand{\Kk}{\mathcal{K}}

\newcommand{\plim}[1]{\displaystyle{\lim_{\stackrel{\longleftarrow}{#1}}}}

\title{On the  Tamagawa number conjecture for Hecke characters}
\author{Francesc Bars \thanks{Work partially supported by BFM2003-06092}}

\begin{document}

\maketitle

\begin{abstract}
In this paper we prove
 the weak $p$-part of the Tamagawa number conjecture in all non-critical cases for the
motives associated to Hecke characters of the form
$\varphi^a\overline{\varphi}^b$ where $\varphi$ is the Hecke
character of a CM elliptic curve $E$ defined over an imaginary
quadratic field $K$,
 under certain restrictions which originate mainly from
 the Iwasawa theory of imaginary quadratic fields.
\end{abstract}





\section{Introduction}
The Tamagawa number conjecture for a variety X over a number field
of
 Bloch and Kato \cite{BK}, or, more
precisely, for a motive $M$ of pure weight $w$ over a number
field, describes the special values of the
 $L$-function in terms of cohomological data (see
for example Kato \cite{Ka1} or Fontaine and Perrin-Riou
\cite{FPR}) and the p-part of the conjecture describes these
values up to units in the ring
$\Z_{(p)}:=\{\frac{a}{b}\in\Q|a,b\in\Z, b\neq0, (b,p)=1\}$.

Recall that the special values of an $L$-function are the leading
coefficient of Taylor expansion at integer
 points. Suppose we have a motive $M$ of weight $w$ such that it's
 $L$-function has meromorphic continuation and satisfies the expected
 functional equation. We say that an integer $m<\frac{w}{2}$ is non-critical
if $L(M,m)=0$ and that it is critical if $L(M,m)\neq0$. We extend
this
 definition to the integers $m>\frac{w}{2}+1$ by saying that $m$ is
critical for $M$ if $w-m+1$ is critical for $\check{M}$ the dual
motive associated to $M$, and non-critical for $M$ if $w-m+1$ is
non-critical for $\check{M}$. The Tamagawa
 number conjecture can be formulated in terms of period
 maps (period integrals nor $p$-adic periods appears) and regulator maps (\cite{Fo},\cite{Ka1}),
 but in the non-critical situation
it can be formulated for almost all non-critical cases (using the
hypothetical functional equation and good compatibilities) without
the period maps, where by period map we mean a map between de Rham
to Betti cohomology and for the $p$-adic period a map between
\'etale cohomology to de Rham cohomology (see \cite[\S 2.3]{Ka1}).

 There are few cases proved in the non-critical situation:  for
  the Riemann zeta function (\cite[\S 6]{BK}), for Dirichlet motives
 (\cite{BG}, \cite{HuKi}), for CM elliptic curves defined over the field of the endomorphism ring
(\cite{K})
 or defined over $\Q$ (\cite[\S7]{BK}, \cite{Ba}).

The weak $p$-part of the Tamagawa number  conjecture for an
elliptic curve $E$
 with CM by $K$ defined over the field of endomorphisms, proved by Kings \cite{K}, is related to
 the weak $p$-part of the Tamagawa number conjecture for the $L$-function of the Hecke character $\varphi$, associated to $E$,
  over the imaginary quadratic field $K$ of class number 1. More precisely, Kings proves
in \cite{K} the conjecture
 for the motive ${h(\overline{\varphi})}(-r)$ with $r\geq0$ which corresponds to the special value (non-critical)
 for the $L$ function associated to $\overline{\varphi}$ at $-r$, where ${h(\overline{\varphi})}$ is the motive associated
 to $\overline{\varphi}$ over $K$ with $K$-coefficients. As a consequence, he obtains the
conjecture for the motive ${h^1(E)}(-r)$. Using the functional
equation for $E$ and good compatibilities one should
 obtain the conjecture for $h^1(E)(r+2)$.
We generalize the methods of Kings to other Hecke characters over
an imaginary quadratic field $K$ in the non-critical
 situation.

We consider the motive associated to the Hecke character
$\varphi^a\overline{\varphi}^b$ with $a,b\geq0$, which has weight
$a+b$. It is known that almost all the non-critical values for
this motive are the integers lower than $min(a,b)$. Our work is
concentrated in this situation, but we remark that there are
results on the Tamagawa number conjecture in the critical
situation (Harrison \cite{Ha}, Guo \cite{Gu}, Kimura \cite{Ki},
Han \cite{Han} and
 in greater generality by Tsuji \cite{Ts}).

 The aim of this paper is to prove the $p$-part of the
 Tamagawa number conjecture (under the formulation in \cite[\S 2.2]{Ka1}) for all the non-critical values
 for the $L$-function of the motive associated Hecke characters
 $\psi_{\theta}$ of an imaginary quadratic field $K$
 with class number 1, which under a fixed embedding corresponds to $\varphi^a\overline{\varphi}^b$.

The main results of this paper are Theorems
\ref{kht:teo},\ref{qht:teo} and Theorems \ref{lnegwTNCk},
\ref{lnegwTNCq}. These results are the weak $p$-part of the
Tamagawa number conjecture for the geometric object associated to
the Hecke characters over $K$ with $a\not\equiv b(modulo
(\#\mathcal{O}_K^*))$ with $K$ or $\Q$-coefficients, under certain
restrictions which originate mainly from Iwasawa theory of
imaginary quadratic fields, where $\mathcal{O}_K^*$ means the
invertible elements of the ring of integers of $K$. See the last
section for numerical examples.

 To obtain these main results, we study in detail the image with respect to
the regulator map of a certain non-trivial submodule of some
$K$-theory group.
 The basic ingredients used in the proof of these results are
Deninger's proof of the Beilinson's conjecture for Hecke
characters in \cite{De}, the specialization of the polylogarithm
sheaf \cite{K} and the Iwasawa main conjecture for imaginary
quadratic fields \cite{Ru3}, as in \cite{K}. This paper need to
deal with negative twists. This problem does not appear in
\cite{K}. For negative twists, we modify Deninger's elements
\cite{De} in order to apply the $p$-adic techniques of \cite{K}.

 The results of this paper generalize the results in \cite[Chapter
 3]{Ba2} which restricts to the Hecke characters $\overline{\varphi}^m$.

\section{The motive associated to Hecke characters}

Let $K$ be an imaginary quadratic field with class number $cl(K)$
equal to 1 and $\Y_K$ be its ring of integers. Let $D_K$ be the
discriminant of $K$. Let $E$ be an elliptic curve over $K$ with CM
by $\Y_K$.
 In this section we describe some pure motives coming from a self
product of the motive $h^1(E)$ and their realizations, and we
prove that the $L$-functions associated to these motives
correspond to Hecke characters. We obtain finally an analog for
these motives of a result of Deuring for CM elliptic curves.

Let $p$ be an odd prime, fixed once and for all, such that $E$ has
good reduction for all primes over $p$. Let $S^{'}$ be the set of
places that divide the conductor of the elliptic curve $\f$ (that
are the same places where $E$ has bad reduction) and the places
that divide $p$.

Let $\underline{\varphi}: I_K\rightarrow K^*$ be the CM character
associated to the elliptic curve $E$ where $I_K$ is the id\`eles
of $K$. Denote by $T_w:=\otimes^w_{\Q}K$ with $w$ a positive
integer. Observe that $T_{w}$ is equal to a product of fields
$\prod_{\theta}T_{\theta}$, where $\theta$ runs through the
$Aut(\C)$-orbits of $\gimel=Hom(T_{w},\C)$, $\theta\subseteq
\gimel$. Let $e_{\theta}$ be the idempotent corresponding to
$T_{\theta}$ of $T_w$.

Define the CM character
$$\underline{\psi}_{\theta}:I_K\rightarrow T_{\theta}^*$$ by
$\underline{\psi}_{\theta}=e_{\theta}\cdot(\otimes^w\underline{\varphi})$,
and denote by $\f_{\theta}$ the conductor of
$\underline{\psi}_{\theta}$. Observe that $\f_{\theta}|\f$ since
$\underline{\psi}_{\theta}$ is a sub-representation of
$\otimes^w\underline{\varphi}$.

Let us fix once and for all an embedding $K\rightarrow\mathbb{C}$
like in the last paragraph on \cite[p.132]{De}.

We have a natural embedding $\yen$
$$K\rightarrow \otimes^w K\rightarrow T_w\rightarrow T_{\theta}$$
where the first map corresponds to the diagonal map.

For any $\vartheta\in\theta$ which its orbit in $\gimel$ is equal
to $\theta$, we have a map $\vartheta:T_w\rightarrow
\vartheta(T_w)\subseteq\mathbb{C}$ inducing an isomorphism
$T_{\theta}=T_w/ker(\vartheta)\rightarrow \vartheta(T_w)$, and
$\vartheta(T_w)$ is the field generated by
$\lambda_1(K)\cdot\ldots\cdot\lambda_w(K)$ where $\lambda_i\in
Hom(K,\C)$, which is isomorphic to $K$.

Define by $\theta_K$ the subset of $\theta$ which contains the
$\vartheta\in\theta$ which orbit is equal to $\theta$ such that
$T_{\theta}=T_w/ker(\vartheta)\rightarrow\vartheta(T_w)\subset\C$
is in $Hom_K(T_{\theta},\C)$ for the natural embedding $\yen$. In
our case $\theta_K$ contains only one element. Let be
$\vartheta=(\lambda_1,\ldots,\lambda_w)\in\gimel$, we set
$a_{\vartheta}=\#\{i|\lambda_i\in Hom_K(K,\C)\}$ and
$b_{\vartheta}=w-a_{\vartheta}$. The infinite type of
$\underline{\psi}_{\theta}$ is defined by
$(a_{\theta},b_{\theta}):=(a_{\vartheta},b_{\vartheta})$ where
$\vartheta$ is the element in $\theta_K$. Observe that $\theta$
only contains two elements, the element of $\theta$ different for
$\vartheta=(\lambda_1,\ldots,\lambda_w)\in \theta_K$ is
$(\overline{\lambda_1},\ldots,\overline{\lambda_w})$ where
$\overline{\lambda_i}$ denotes the composition of $\lambda_i$ with
the complex conjugation.

Consider the category of Chow motives $\mathcal{M}(K)$ over $K$
with morphisms induced by graded correspondences in Chow theory.
We have then a natural covariant functor $h$ from the category of
smooth and projective varieties over $K$ to $\mathcal{M}(K)$. Then
the motive $h(E)$ of an elliptic curve $E$ over $K$ has a
decomposition with respect to the zero section $h(E)\cong h^0E
\oplus h^1E \oplus h^2E$, where $h^0E=h(Spec(K))$ and
$h^2E=h(Spec(K))(-1)$. We can also consider the category
$\mathcal{M}_{\Q}(K)$ which consist of the same objects but
tensoring by $\Q$ the group of morphisms.

The motive $h^1E$ has multiplication by $\Y_K$. Consider then the
motive $\otimes^w h^1E,$ has multiplication by
$\Y_{w}:=\otimes^w_{\Z}\Y_K$. Then $\otimes^w h^1E_{\Q}$ has
multiplication by $T_w$. Notice that $e_{\theta}$ is not integral
in general for $w>1$, but is contained in $\Y_K[1/D_K]$. Let's
denote by
$$M_{\theta}:=e_{\theta}(\otimes^wh^1(E)_{\Q}\otimes_{\Y_K}\Y_K[1/D_K]),$$
considered as an motive with coefficients in $\Y_K[1/D_K]$, and by
$M_{\theta\Q}$ its image in $\mathcal{M}_{\Q}(K)$. As
$\lambda_i\in\{\lambda,\overline{\lambda}\}$ where $\lambda$ is
the fixed embedding of $K$ in $\C$, and $\lambda_i(\Y_K)=\Y_K$, we
have then that $e_{\theta}(\otimes
h^1(E)\otimes_{\Y_K}\Y_K[1/D_K])\in\mathcal{M}(K)$ has
multiplication by
$\Y_{\theta}:=e_{\theta}(\Y_w\otimes_{\Y_K}\Y_K[1/D_K])$, and
$M_{\theta\Q}$ has multiplication by $T_{\theta}\cong K$.


Our objective in this section is to study the $p$-adic and Betti
realizations of this motive $M_{\theta}$, called Hecke motive, and
to determine its $L$ function.

There are at least three equivalent notions of a Hecke character,
see \cite[p.48]{G}. One is the notion of CM-character used above,
\cite[p.48, definition 2]{G}. For Hecke $L$-functions and the
Galois group action on the $p$-adic realization associated to the
Hecke motive, we use the notion of a character which is trivial on
$K^*$ and with image in some id\`ele group,
$\psi_{\theta}:I_K/K^*\rightarrow I_{T_{\theta}}$ \cite[p.48,
definition 3]{G}. The associated complex Hecke character, in order
to define the Hecke $L$-function, is constructed from
$\psi_{\theta}$ by taking the archimedian places of
$I_{T_{\theta}}$ which correspond to the fixed immersion of $K$ in
$\C$ in our situation, which we also call $\psi_{\theta}$. The
character constructed from $\psi_{\theta}$ by taking the
components of the places of $I_{T_{\theta}}$ above $p$ is called
$\psi_{\theta,p}$ and is related with the Galois action on the
$p$-adic realization associated to the motive.
The character $\psi_{\theta,p}$ factors through $Gal(K^{ab}/K)$
and has image in $(T_{\theta}\otimes\Z_p)^*$. We will use the term
Hecke character when we want to consider this second notion from
now on. The third notion \cite[definition 1,p.48]{G}, corresponds
to certain map
$\widetilde{\psi}_{\theta}:\mathcal{I}_{\f_{\theta}}\rightarrow
T_{\theta}^*$, where $\mathcal{I}_{\f_{\theta}}$ is the ideal
classes of $K$ prime to $\f_{\theta}$. So, if $\p$ is a prime
ideal of $\mathcal{O}_K$ prime to $\f_{\theta}$, we mean for
$\psi_{\theta}(\p)$ or $\underline{\psi}_{\theta}(\p)$ the value
of the Hecke character or CM character at the id\`ele which has an
uniformizer $\pi$ at the place $\p$ and 1 in the other places. We
have
$\psi_{\theta}(\p)=\underline{\psi}_{\theta}(\p)=\widetilde{\psi}_{\theta}(\p)$
(see \cite[p.49-50]{G}).

The $p$-adic realization of the motive $M_{\theta\Q}(w)$ is, by
definition, $H^w_{et}(M_{\theta\Q}\times_{K}\overline{K},\Q_p(w))$
and we denote it by $M_{\theta\Q_p}(w)$.

\begin{lem} Let $p$ be a prime such that $p\nmid D_K$.
The integral $p$-adic realization of $M_{\theta}(w)$,
$H^w_{et}(M_{\theta}\times_{K}\overline{K},\Z_p(w))\otimes_{\Y_K}\Y_K[1/D_K]$,
is isomorphic
 to $$e_{\theta}(\otimes^wT_pE)$$ as  free
$e_{\theta}(\otimes^w\Y_K[1/D_K]\otimes\Z_p)$-modules of rank 1,
with $Gal(\overline{K}/K)$-action on $e_{\theta}(\otimes^wT_pE)$
given by the Hecke character $\overline{\psi}_{\theta,p}$.
\end{lem}
\begin{proof} Observe first that $T_pE$ is isomorphic as Galois modules to
$H^1_{et}(h^1(E)\times_K\overline{K},\Z_p(1))=Hom(T_pE,\Z_p(1))$
by the use of Weil pairing but the isomorphism change the action
of $\mathcal{O}_K$ to its conjugate.

The claim that $e_{\theta}(\otimes^wT_pE)$ is a free module of
rank 1 follows
 because $T_pE$ is a free $\Y_K\otimes\Z_p\cong \Y_K[1/D_K]\otimes\Z_p$-module of rank 1 and then
$e_{\theta}\cdot (T_pE\otimes\ldots\otimes T_pE)$ is a free
$e_{\theta}\cdot(\otimes^w(\Y_K[1/D_K]\otimes\Z_p))$-module of
rank one.

Now, consider the natural action of $G_K:=Gal(\overline{K}/K)$ on
$H^1_{et}(h^1(E)\times_K\overline{K},\Z_p(1))$. Since $G_K$ acts
on the Tate module by the Hecke character
$\varphi_p:G_K\rightarrow(\Y_K\otimes\Z_p)^*$, so it acts on
$H^1_{et}(h^1(E)\times_K\overline{K},\Z_p(1))$ by
$\overline{\varphi}_p$. Using
$$H^w((h^1(E)\times_K\overline{K})^w,\Z_p(w))=
H^1(h^1(E)\times_{K}\overline{K},\Z_p(1))^{\otimes w}$$ and taking
our idempotent, the action is given by
$e_{\theta}(\otimes^w\overline{\varphi}_p)=\overline{\psi}_{\theta,p}$.
\end{proof}
Twisting by $l+1$ we get for $p\nmid D_K$ that the integral
$p$-adic realization for $M_{\theta}(w+l+1)$ is isomorphic to
$e_{\theta}(\otimes^w T_pE)(l+1)$ with $G_K$-action on
$e_{\theta}(\otimes^w T_pE)(l+1)$ given by
$\overline{\psi}_{\theta,p}$ multiplied by the $l+1$-th power of
the $p$-adic cyclotomic character.

The Betti realization for the motive $M_{\theta\Q}(w+l)$ named
$H^w_B(M_{\theta\C},\Q(w+l))$ is isomorphic to
$e_{\theta}(\otimes^wH^1_B(E(\C),\Q(1)))(l)$, we remember that we
fixed once and for all an immersion for $K\subseteq\C$. $E(\C)$ is
the set of closed points with the analytic topology. We have
$$\otimes^w_{\Z}(H^1_B(E\times_K
\C),\Z(1))\otimes_{\Y_K}\Y_K[1/D_K]$$ a $\otimes^w
\Y_K[1/D_K]$-module of rank 1 and taking the idempotent
$e_{\theta}$ we obtain $$e_{\theta}(\otimes^w
H^1_B(E\times_K\C),\Z(1))\otimes_{\Y_K}\Y_K[1/D_K])(l),$$ a
$\Y_{\theta}$-module of rank 1 which is the submodule of
$H^w_B(M_{\theta\C},\Q(w+l))$ corresponding to
$H^w_B(M_{\theta\C},\Z(w+l)))\otimes_{\Y_K}\Y_K[1/D_K]$.

Now, we are going to study the $L$-function that corresponds to
 the $p$-adic representation $M_{\theta\Q_p}=H^w_{et}(M_{\theta\Q}\times_K\overline{K},\Q_p) $
of $M_{\theta\Q}$.

The Tamagawa number conjecture describes conjecturally special
values of the $L$-function for the motive and this $L$-function
involves the product of all Euler factors (one for every
non-arquimedian place), but for the $p$-part of this conjecture,
Kato reformulates the conjecture in terms of the partial
$L$-function avoiding a non-vanishing finite set of Euler factors,
more concretely the finite set contains the Euler factors coming
from the primes above $p$ and the primes where the motive has bad
reduction (see \cite[Proposition 7.8]{Ka2}, or \cite[Chapter
1]{Ba2} for an overview, and see remark \ref{rem26} for these
no-vanishing in our setting).

Let $S$ be the set of places of $K$ that divide $\f_{\theta}$ or
$p$. Define as usual
$$L_S(M_{\theta\Q},s):=
\prod_{\mathfrak{l}\notin S}
det_{\Q_p}(1-Frob_{\mathfrak{l}}N\mathfrak{l}^{-s}|
M_{\theta\Q_p}^{I_{\mathfrak{l}}})^{-1}$$ where
$Frob_{\mathfrak{l}}$ means the geometric Frobenius,
$I_{\mathfrak{l}}$ the inertia group at $\mathfrak{l}$ and
$N\mathfrak{l}$ the norm $N_{K/\Q}\mathfrak{l}$.

Our goal is to compute this determinant and to relate it to the
local factors of the $L$-function of the Hecke character
$\psi_{\theta}$ that is defined by
$$L_S(\psi_{\theta},s):=\prod_{\mathfrak{l}\notin S}
(1-\frac{\psi_{\theta}(\s)}{N\s^s})^{-1}.$$

Recall that the operation of the decomposition group $D_{\p}$ for
$\p\nmid p$ on $H^1_{et}(h^1(E)\times_K\overline{K},\Q_p)$ is
given by the operation of $\varphi^{-1}|_{K_{\p}^*}$, and hence
$D_{\p}$ operates on $M_{\theta\Q_p}$ via $\psi_{\theta}^{-1}$. On
one hand, the inertia group $I_{\p}$ acts non-trivially if and
only if $\p$ divides the conductor $\f_{\theta}$. On the other
hand, for $\p\nmid\f_{\theta}$, the geometric Frobenius $Fr_{\p}$
at $\p$ acts via $\psi_{\theta}(\p)$. We obtain hence the
following result.

\begin{lem}[Deninger, prop. 1.3.2.a \cite{De}] \label{den:teo} Let $\s$ be a finite prime of $K$ with $\s\nmid p\f_{theta}$,
where $\mathfrak{f}_{\theta}$ is the conductor of the Hecke
character $\psi_{\theta}$. Then
$$det_{T_{\theta}\otimes\Q_p}(1-Fr_{\s}N\s^{-s}|M_{\theta\Q_p}^{I_{\s}})
=(1-\psi_{\theta}(\s)N\s^{-s}).$$
\end{lem}

We impose some restrictions for our motive $M_{\theta}(w+l+1)$
once and for all. We suppose $-w-2l\le-3$. Remember that, with the
restriction that $E$ is defined over $K$, we have $\#|\theta|=2$,
and, in particular, we have $T_{\theta}\cong K$ and for $p\nmid
D_K$ we have $\Y_{\theta}\otimes\Z_p\cong\Y_K\otimes\Z_p$.

The $L$-function for $M_{\theta\Q}$ can be described by using
lemma \ref{den:teo} and by taking the norm map.
\begin{lem}\label{spl:teo}
Let $\s$ a prime of $K$ such that $\s\nmid\mathfrak{f}_{\theta}$
and it is prime to $p$. We have then the following equality
$$det_{\Q_p}(1-Fr_{\s}N\s^{-s}|M_{\theta\Q_p}^{I_{\s}})=
(1-\psi_{\theta}(\s)N\s^{-s})(1-\overline{\psi_{\theta}}(\s)N\s^{-s}).$$
\end{lem}

As a corollary we obtain a generalization of a result of Deuring.

\begin{thm}\label{hdeu:teo} Let $S$ be the set of the
primes on $K$ dividing $\f_{\theta}$ and primes dividing $p$.
Then:
$$L_S(M_{\theta\Q},s)=L_S(\psi_{\theta},s)L_S(\overline{\psi}_{\theta},s).$$
\end{thm}

\begin{remark}\label{rem26} The $p$-adic realization
$V_p:=M_{\theta\Q_p}(w+l+1)$ satisfies that the local Euler
factors
$$P_{\mathfrak{p}}(V_p^*(1),0)=P_{\mathfrak{p}}(\overline{\psi}_{\theta},-l)$$
are different from 0 for all $\mathfrak{p}\in S$ where $V_p^*$ is
the dual Galois module of $V_p$. Hence, it satisfies the
hypothesis of \cite[conjecture 2.2.7]{Ka1}.

To show this fact, suppose first that
$\mathfrak{p}|\mathfrak{f}_{\theta}$. Then the inertia group acts
non-trivially on $V_{p}$, which is a one dimensional
$\Y_{\theta}\otimes\Q$-module, and hence
$$L_{\mathfrak{p}}(\overline{\psi}_{\theta},s)=1$$
for all $\mathfrak{p}|\mathfrak{f}_{\theta}$, and in particular
for $s=-l$.

If $\mathfrak{p}$ divides $p$, then the result follows from the
fact that any proper smooth variety with good reduction at
$\mathfrak{p}$ satisfies it for weight reasons, and in particular:
$$det_{\Q_p}(1-Fr_{\mathfrak{p}}N\mathfrak{p}^{l}|H^w_{et}((E\times_K\overline{K})^w,\Q_{l'}))\neq0,$$
with $l'\neq p$, and therefore, since the different idempotents
$e_{\theta}$ give a direct summand of the cohomology group
$H^w((E\times_K\overline{K})^w,\Q_{l'})$,
$$L_{\mathfrak{p}}(\overline{\psi}_{\theta},-l)\neq0.$$
\end{remark}

The motivic cohomology group
$H_{\Mm}^{w+1}(M_{\theta\Q},\Q(w+l+1))$ is the $K$-theory group
$K_{2(w+l)-w+1}(M_{\theta\Q})^{(w+w+1)}\otimes\Q$ where the
$K$-groups are the Quillen $K$-groups and the superscript denotes
the Adam's filtration on them.

We suppose recall that $w-2(w+l+1)\leq-3$. We have a Beilinson
regulator map,
$$r_{\Dd}:H_{\Mm}^{w+1}(M_{\theta\Q},\Q(w+l+1))\rightarrow H_B^w(M_{\theta\C},\Q(w+l))\otimes_{\Q}\R,$$
and the Soul\'e regulator map:
$$r_p: H_{\Mm}^{w+1}(M_{\theta\Q},\Q(w+l+1))\rightarrow
H^1_{\acute{e}t}(\Y_K[1/S],M_{\theta\Q_p}(w+l+1)).$$

The $p$-part of the Tamagawa number conjecture claims in
particular that $r_{\Dd}\otimes_{\Q}\R$ and $r_p\otimes_{\Q}\Q_p$
are isomorphisms. Deninger in \cite{De}\cite{Ded} constructs a
$\Q$-subspace $H_{\Mm}^{constr}$ of
$H_{\Mm}^{w+1}(M_{\theta\Q},\Q(w+l+1))$ such that
$r_{\Dd}\otimes\R$ restricted to $H_{\Mm}^{constr}\otimes_{\Q}\R$
is an isomorphism. The term weak in the formulation of the
$p$-part of the Tamagawa number conjecture means that this
conjecture is proved using the space $H_{\Mm}^{constr}$ instead of
$H_{\Mm}^{w+1}(M_{\theta\Q},\Q(w+l+1))$.

\section{The Beilinson conjecture for Hecke characters}
In this section we review the work on the Beilinson conjecture for
the motive $M_{\theta\Q}(w+l+1)$ done by Deninger in \cite{De},
under the language of the Tamagawa number conjecture.

\begin{thm}[Deninger, Theorem 1.4.1 \cite{De}]\label{nul:teo} Let $w=a_{\theta}+b_{\theta}\ge1$. Consider an integer $l$ such that
$$-l\leq Min(a_{\theta},b_{\theta})\ \ \ if\ a_{\theta}\neq b_{\theta}$$
$$-l<a_{\theta}=b_{\theta}=w/2\ \ \ otherwise.$$
Then the L-series $L(\overline{\psi}_{\theta},s)$ has a zero of
order 1 at $s=-l$, (i.e.
$ord_{s=-l}L(\overline{\psi}_{\theta},s)=1$).

Moreover, there exist an element $\xi_{\theta}$ in
$H_{\Mm}^{w+1}(M_{\theta\Q},\Q(w+l+1))$, such that
$$r_{\Dd}(\xi_{\theta})=
\lim_{s\rightarrow
-l}\frac{L(\overline{\psi_{\theta}},s)}{s+l}\eta_{\theta} \ mod\
T_{\theta}^*$$ in the free rank one $T_{\theta}\otimes\R$-module
$H_B^w(M_{\theta\C},\R(w+l))$, where $\eta_{\theta}$ is a
$T_{\theta}$-generator of $H_B^w(M_{\theta\C},\Q(w+l))$.
\end{thm}

Let's recall the construction of $\xi_{\theta}$, following the
results of Deninger. We suppose once for all that $l\ge0$.

Fix an algebraic differential form $\omega\in
H^0(E,\Omega_{E/K})$. Since we have complex multiplication, we can
write the period lattice as $\Gamma=\Omega\Y_K$, where
$\Omega\in\C^*$ is the complex period. Fix an element $\gamma$ in
$H_1(E(\C),\Z)$ such that it is an $\Y_K$-generator, and satisfies
$$\Omega=\int_{\gamma}\omega.$$
By Poincar\'e duality, we have that $\gamma$ corresponds to
$\eta_{\gamma}$, an $\Y_K$-generator for $H^1(E(\C),\Z(1))$. Thus
$\eta_{\gamma}\otimes_{\Y_K}\Y_K[1/D_K]$ is an
$\Y_K[1/D_K]$-generator for the module
$H^1(E(\C),\Z(1))\otimes_{\Y_K}\Y_K[1/D_K]$ which by abuse of
notation we call also $\eta_{\gamma}$. Consider now the
$\Y_{\theta}$-generator
$$\eta_{\theta}:=(2\pi i)^l e_{\theta}(\otimes^w\eta_{\gamma})$$
of $H_B^w(M_{\theta\C},\Z(w+l))\otimes_{\Y_K}\Y_K[1/D_K]$.

To construct $\xi_{\theta}$, we will define a divisor on the
torsion points of the elliptic curve; its image by the composition
of the Eisenstein map $\mathcal{E}_{\Mm}$ (\cite[\S 8]{Ded}) with
the Deninger projector map $\mathcal{K}_{\Mm}$ (\cite[(2.8)]{De})
will define our $\xi_{\theta}$.

Remember that $\f_{\theta}$ is the conductor of the Hecke
character $\psi_{\theta}$ associated with $M_{\theta}$, and denote
by $f$ a generator of $\f_{\theta}$ (it exists since $cl(K)=1$).
We have that
$$\Omega f^{-1}\in \f_{\theta}^{-1}\Gamma$$
and that $(\Omega f^{-1})$ gives a divisor in
$\Z[E[\f_{\theta}]\setminus0]$ defined over $K(E[\f_{\theta}])$.
Since $\mathfrak{f}$ is the conductor of $\psi$ and
$\f_{\theta}|\f$, the divisor $(\Omega f^{-1})$ is defined also
over $K(E[\mathfrak{f}])$. We will define our divisor as
 $$\beta_{\theta}:=N_{K(E[\f])/K}((\Omega f^{-1})).$$
 Denote by $\rho_{\theta}$ a finite id\`ele such that
 $(\rho_{\theta})=\f_{\theta}$ and
 $v_{\p}(f^{-1}-\rho_{\p}^{-1})\geq 0$ for all $\p\mid
 \f_{\theta}$.

If $a_{\theta}\not\equiv b_{\theta}$ mod $|\Y_K^*|$, we obtain
that (\cite[p.142,(2.11)]{De})
$$r_{\Dd}(\mathcal{K}_{\Mm}\mathcal{E}_{\Mm}(\beta_{\theta}))=
(-1)^{l-1}\frac{2^{l-1}N_{K/\Q}\f_{\theta}^{w+2l}\psi_{\theta}(\rho_{\theta})}
{(2l+w)!N_{K/\Q}(\f_{\theta})^{l+w}}\frac{\Phi(\f)}{\Phi(\f_{\theta})}
L'(\overline{\psi_{\theta}},-l)\eta_{\theta},$$ where
$\Phi(\mathfrak{m}):=|(\Y_K/\mathfrak{m})^*|$ for any ideal
$\mathfrak{m}$ of $\Y_K$.

This is an analog for $M_{\theta}(w+l+1)$ of \cite[thm.1.2.2]{K},
which corresponds to the case $h^1(E)(1+l+1)$.

\begin{thm}[Deninger, \S 2 \cite{De}]\label{denig:teo} Suppose that $a_{\theta}\not\equiv b_{\theta}$ mod $(\#\Y_K^*)$ and
that $a_{\theta}, b_{\theta}, l$ satisfy the hypothesis of the
theorem \ref{nul:teo} with $l\ge0$. Define, by using the previous
notation, $$\xi_{\theta,l}:=$$
$$(-1)^{l-1}\frac{(2l+w)!L_p(\overline{\psi_{\theta}},-l)^{-1}
\Phi(\f_{\theta})}{2^{l-1}N_{K/\Q}\f_{\theta}^{l}\psi_{\theta}(\rho_{\theta})\Phi(\f)}\mathcal{K}_{\Mm}
\circ\mathcal{E}_{\Mm}^{2l+w}(\beta_{\theta})\in
H^{w+1}_{\Mm}(M_{\theta\Q},\Q(w+l+1)),$$ where
$L_p(\overline{\psi_{\theta}},s)$ is the product of the Euler
factors for the primes of $K$ above $p$.

 Then
$$r_{\Dd}(\xi_{\theta,l})=L_S^*(\overline{\psi_{\theta}},-l) e_{\theta}(\otimes^w\eta_{\gamma}),$$
where $S$ are the primes of $K$ that divide $\f_{\theta}p$. Here
$L_S^*(\overline{\psi_{\theta}},-l)=\lim\limits_{s+l\rightarrow 0}
L_S(\overline{\psi_{\theta}},s)/(s+l)$.
\end{thm}

\begin{definition} For $a_{\theta}\not\equiv b_{\theta}$ mod$(\# \Y_K^*)$ we define
$$\Rr_{\theta}:=\xi_{\theta,l}\Y_K.$$
\end{definition}

\begin{remark}
Theorem 1.4.1 \cite{De} is more general because it includes the
situation $a_{\theta}\equiv b_{\theta} \ mod(\#\Y_K^*)$. But in
this situation, Deninger defines a divisor
$\tilde{\beta_{\theta}}$ instead of $\beta_{\theta}$ which is not
a norm of a positive divisor and moreover it contains the zero
point of $E$. Thus the techniques for constructing an Euler system
of \S 5 can not be applied in this case (see for example Theorem
\ref{poleli}).
\end{remark}

As a consequence of Theorem \ref{denig:teo}, we have that our
submodule $\Rr_{\theta}$ verifies some integral version of the
Beilinson conjecture for the motive $M_{\theta}(w+l+1)$.

\begin{thm} \label{yok:teo}
The $\Y_{K}$-submodule $\Rr_{\theta}$ of
$H^{w+1}_{\Mm}(M_{\theta\Q},\Q(w+l+1))$ satisfies that

 $$det_{\Y_{K}[1/D_K]}(r_{\Dd}(\Rr_{\theta}\otimes_{\Y_K}\Y_K[1/D_K])) = $$
 $$L_S^*(\overline{\psi_{\theta}},-l)
det_{\Y_{K}[1/D_K]}(H^w_B(M_{\theta},\Z(w+l))\otimes_{\Y_K}\Y_K[1/D_K])$$

in
$det_{\Y_{K}[1/D_K]\otimes\R}((H^w_B(M_{\theta\C},\Z(w+l))\otimes_{\Y_K}\Y_K[1/D_K])\otimes\R)$.

\end{thm}
\begin{proof} Observing that $\eta_{\theta}$ is
a $\Y_{\theta}$-base for the free $\Y_{\theta}$-module
$$H^w_B(M_{\theta\Q}\times_K\C,\Z(w+l))\otimes_{\Y_K}\Y_K[1/D_K]$$
of rank one, the result follows.
\end{proof}
\begin{cor} \label{yo:teo}
The submodule $\Rr_{\theta}$ defined above satisfies the Beilinson
conjecture inside the $p$-part of the Tamagawa number conjecture
for $p\nmid D_K$,
 that is $\Rr_{\theta}$ satisfies the following conditions:
\begin{enumerate}
\item The map $r_{\Dd}\otimes\R$ is a isomorphism when  restricted
to $\Rr_{\theta}\otimes\R$. \item
$dim_{\Q}(H^w_B(M_{\theta\C},\Q(w+l)))=ord_{s=-l}L_S(M_{\theta\Q},s)=2$.
\item We have the following equality
$$r_{\Dd}(det_{\Z[1/D_K]}(\Rr_{\theta}\otimes_{\Y_K}\Y_K[1/D_K]))=$$
$$L_S^*(M_{\theta\Q},-l)
det_{\Z[1/D_K]}(H^w_B(M_{\theta},\Z(w+l))\otimes_{\Y_K}\Y_K[1/D])$$
where $L_S^*(M_{\theta\Q},-l)$ means $\lim\limits_{s\rightarrow
-l}L_S^*(M_{\theta\Q},s)/(s+l)^2$ (this makes sense by using
theorem \ref{hdeu:teo} and theorem \ref{nul:teo}).
\end{enumerate}
\end{cor}
\begin{proof} The first and the second conditions are clear for the dimensions
of the spaces involved in the Deligne regulator map, and the
theorem \ref{yok:teo}. The third condition comes from the previous
theorem using the fact that, if we multiply an
$\Y_{\theta}$-module with an element
$L_S^*(\overline{\psi_{\theta}},-l)$ in $\Y_{\theta}\otimes\R$,
the determinant is multiplied by the norm
$$N_{\Y_{\theta}\otimes\R/\R}(L_S^*(\overline{\psi_{\theta}},-l))=
L_S^*(\overline{\psi_{\theta}},-l)\overline{L_S^*(\overline{\psi_{\theta}},-l)}=
L_S^*(\overline{\psi_{\theta}},-l)L_S^*(\psi_{\theta},-l).$$ Using
theorem \ref{hdeu:teo}, we obtain that this is equal to
$L_S^*(M_{\theta\Q},-l)$.
\end{proof}

\section{Iwasawa theory}

We suppose once and for all that $p\nmid \# \Y_K^*$ and $p\nmid
N_{K/\Q}(\f)$, (in particular $p\nmid D_K$).

To simplify the notation, we will denote in the following by
$$M_{\theta\Z_p}(w+l+1)=e_{\theta}(\otimes^w (H^1_{et}(E\times_K\overline{K},\Z_p(1)))(l+1)$$
the $p$-adic lattice for the $p$-adic realization of
$M_{\theta}(w+l+1)$.

Let $K_n:=K(E[p^{n+1}])$ be the field of definition of the
$p^{n+1}$-torsion points of $E$, $\Y_n$ its ring of integers and
let $K_{\infty}:=\lim\limits_{\rightarrow}K_n$ be its direct
limit. Denote by $\Y_n$ the ring of integers of $K_n$
(respectively $\Y_{\infty}$). We know that $\Delta:=Gal(K_0/K)$
has order prime to $p$ and $\Gamma:=Gal(K_{\infty}/K_0)$ is
isomorphic to $\Z_p^2$.

Let $\mathcal{G}$ be the Galois group $Gal(K_{\infty}/K)$; then
$\mathcal{G}\cong\Delta\times\Gamma$.

We use now the notations on Iwasawa theory for imaginary quadratic
fields used in \cite[\S 2.1]{K} but with a different definition of
elliptic units.

Let us define the elliptic units $\Cc_{n,\f_{\theta}}$ in $K_n$
which are more convenient for us. \\

For every ideal $\mathfrak{a}$ of $K$ prime to $6$ we can define a
theta function
$$\theta_{\mathfrak{a}}: E\setminus ker([\mathfrak{a}])
\longrightarrow \mathbb C$$ which has divisor
$N(\mathfrak{a})(e)-ker([\mathfrak{a}])$ (for the precise
definition see \cite[Theorem 4.2.2]{K}). The function
$\theta_{\mathfrak{a}}(z)$ is in fact a 12-th root of the function
defined in \cite[II.2.4]{dS}. Let $\mathfrak{g}$ be a fixed ideal
of $\Y_K$ such that $\Y_K^*\rightarrow(\Y_K/\mathfrak{g})^*$ is
injective, and suppose that $\mathfrak{g}$ divides the conductor
$\f$ of the elliptic curve $E$. Let's denote by $t_{\mathfrak{g}}$
a generator for the $E[\mathfrak{g}]$-torsion points as
$\Y_K$-module, and let $\mathfrak{a}$ be an ideal prime to
$6\mathfrak{g}$.
\begin{definition}\label{eliuni}Let $C_{n,\mathfrak{g}}$ be the subgroup of units generated over $\Z[Gal(K_n/K)]$
by $$\prod_{\sigma\in
Gal(K(\mathfrak{g})/K)}\theta_{\mathfrak{a}}(t_{\mathfrak{g}}^{\sigma}+h_n),$$
where $\mathfrak{a}$ runs through all ideals prime to $6p\f$,
$K(\mathfrak{g})$ is the ray class field defined by $\mathfrak{g}$
and $h_n$ is a primitive $p^{n+1}$-torsion point (i.e. a generator
of the $p^{n+1}$-torsion points of $E$ as $\Y_K$-module). Define
the group of elliptic units of $K_n$ as
$$\Cc_{n,\mathfrak{g}}:=\mu_{\infty}(K_n)C_{n,\mathfrak{g}},$$
where $\mu_{\infty}(K_n)$ denotes the roots of unity in $K_n$.
\end{definition}

Denote by $\overline{\Cc}_{n,\mathfrak{g}}$ the closure in the
local units $\mathcal{U}_n^{\mathfrak{p}}$ where
$\mathcal{U}_n^{\mathfrak{p}}$ is the group of local units of
$K_n\otimes_{K}K_{\mathfrak{p}}$ which are congruent to 1 modulo
the primes above $\mathfrak{p}$ where $\mathfrak{p}$ is a prime of
$K$ above $p$. Define $\overline{\Cc}_{\infty,\mathfrak{g}}:=
\lim\limits_{\leftarrow}\overline{\Cc}_{n,\mathfrak{g}},$ and
$\mathcal{U}_{\infty}^{\mathfrak{p}}:=\lim\limits_{\leftarrow}\mathcal{U}_n^{\mathfrak{p}}$
where the limit is taken with respect to the norm maps. Define
also $\mathcal{U}_{\infty}$ by
$\mathcal{U}_{\infty}^{\mathfrak{p}}\times\mathcal{U}_{\infty}^{\mathfrak{p}^*}$
if $p=\mathfrak{p}\mathfrak{p}^*$ splits, and if $p$ inert or
ramified by $\mathcal{U}_{\infty}^{\mathfrak{p}}$. Let
$\mathcal{Y}_n$ be the $p$-adic completion of $(K_n\otimes\Q_p)^*$
and $\mathcal{Y}_{\infty}:=\lim\limits_{\leftarrow}\mathcal{Y}_n$.


Let us once and for all to specialize the elliptic units to
$\mathfrak{g}=\f_{\theta}$.

 Recall that $S$ denotes the set of primes
 of $K$ which divide $\f_{\theta}$ or $p$, and that $S'$ denotes the set
 of primes of $K$ which divide $p$ or the conductor
$\f$ of the
 elliptic curve $E$. Denote $\Y_S:=\Y_K[1/S]$ and $\Y_p:=\Y_K\otimes\Z_p$.

We are going to define a map in the spirit of Soul\'e:
$$(Soul)_p:\overline{\mathcal{C}}_{\infty,\f_{\theta}}\otimes_{\Z_p}
M_{\theta\Z_p}(w+l)\rightarrow H^1(\Y_S,M_{\theta\Z_p}(w+l+1)),$$
observe that $M_{\theta\Z_p}(w+l)$ is unramified outside $S$, thus
$H^1(\Y_S,M_{\theta\Z_p}(w+l+1))$ makes sense.

Using the definition of $M_{\theta\Z_p}(w)(l+1)$, we have that
$$H^1(\Y_S,M_{\theta\Z_p}(w+l+1))=
\lim\limits_{\leftarrow}H^1(\Y_S,(e_{\theta}\otimes^wE[p^{r+1}])(l+1)).$$
Define $(Soul)_p$ in the following way. Given $(\theta_r)_r$ a
norm compatible system of elliptic units and an element $(t_r)_r$
of $\lim\limits_{\leftarrow}(e_{\theta}(\otimes^wE[p^{r+1}]))(l)$,
we define
$$(Soul)_p((\theta_r\otimes t_r)_r):=(N_{K_r/K}(\theta_r\otimes t_r))_r .$$
It is well defined: $\theta_r\otimes t_r$ is an element in
$$\Y_{r,S}^*/(\Y_{r,S}^*)^{p^{r+1}}\otimes
(e_{\theta}(\otimes^wE[p^{r+1}]))(l)\subset
H^1(\Y_{r,S},(e_{\theta} (\otimes^wE[p^{r+1}]))(l+1))$$ where
$\Y_{r,S}$ is $\Y_r[1/S]$ the ring of integers of $K_r$ inverting
the primes above $S$, $N_{K_r/K}$ denotes the norm map in
cohomology and by Soul\'e's Lemma 1.4 \cite{So} one gets an
element in $H^1(\Y_S,M_{\theta\Z_p}(w+l+1))$. The map $(Soul)_p$
factors thought the coinvariants, denoted by
$(\overline{\Cc}_{\infty,\f_{\theta}}\otimes
M_{\theta\Z_p}(w+l))_{\mathcal{G}}$.
\begin{definition} The Soul\'e elliptic elements are the elements in the image of the map
$$(Soul)_p:(\overline{\Cc}_{\infty,\f_{\theta}}\otimes M_{\theta\Z_p}(w+l))_{\mathcal{G}}
\rightarrow H^1(\Y_S,M_{\theta\Z_p}(w+l+1))$$ where
$\mathcal{G}=Gal(K(E[p^{\infty}])/K)$.
\end{definition}

We consider in the following the representation $\chi$ of the
group $\Delta$ given by the action of $\Delta$ in
$Hom_{\Y_p}(M_{\theta\Z_p}(w+l),\Y_p)$.


We are only able to apply the techniques on Iwasawa theory of
\cite{K} for certain representations that we call good
representation.
\begin{definition}\label{goo:teo}\label{4.3}
We say that such representation $\chi$ of the group $\Delta$ is a
good representation if it satisfies two conditions in Iwasawa
theory about isomorphism between some concrete Iwasawa modules:
(A) the Iwasawa main conjecture of Rubin \cite[Theorem 2.1.3]{K}
but replacing the elliptic units module there with the elliptic
module units $\overline{\mathcal{C}}_{\infty,\f_{\theta}}^{\chi}$
and (B) from the inclusion
$\mathcal{U}_{\infty}\subseteq\mathcal{Y}_{\infty}$ we get that
$\mathcal{U}_{\infty}^{\chi}\cong \mathcal{Y}_{\infty}^{\chi}$ as
Iwasawa modules for the Iwasawa ring
$\lim\limits_{\leftarrow}\Z_p[[Gal(K_n/K)]]^{\chi}$.
\end{definition}
We observe that the elliptic units $\Cc_{\infty,\f}$, which are
the ones that appear in \cite{Ru4} and \cite{K}, satisfies the
theorem of Iwasawa main conjecture of \cite{Ru3} for any
$\Delta$-representation under the hypothesis of the theorem in
\cite{Ru3} (personal communication of Rubin).

When $S=S'$ we have $\Cc_{\infty,\f_{\theta}}=\Cc_{\infty,\f}$
\cite[Proposition II.2.5]{dS}, therefore the Iwasawa main
conjecture (condition (A)) is true from Rubin's theorem
\cite{Ru3}\cite{Ru4} for $p$ splits and for $p$ inert when $\chi$
is non trivial on the decomposition group $\Delta_{\mathfrak{p}}$
of $\mathfrak{p}$ in $\Delta$.

Condition (B) is always true if $p$ splits \cite[Lemma 2.1.6]{K}.
 If $p$ is inert or ramified, the representation $\chi$ satisfies condition (B)
 if $\Z_p[\Delta/\Delta_p]^{\chi}=0$ (see \cite[lemma 2.1.6]{K}), moreover because $p$ is a prime over
 which $E$ has good reduction we have $\Delta_p=\Delta$ \cite[Lemma
2.2.9]{K}, thus condition (B) is true if $\chi$ is non-trivial.

\begin{remark} We guess that the Iwasawa main conjecture \cite[Theorem 2.1.3]{K} is also
true without the condition $S=S'$ for our elliptic units and our
character $\chi$, so such a character should be good if it just
verify condition (B). In the next section we prove that
$\overline{\Cc}_{\infty,\f_{\theta}}^{\chi}$ is a rank 1 Iwasawa
module and we construct an Euler system there. Using this and
\cite[Lemma III.1.10]{dS}, it should be possible to prove the
Iwasawa main conjecture  for $\chi$ using the techniques in
\cite{Ru4}.
\end{remark}

\begin{lem}\label{4.5} Suppose that $\psi_{\theta}$ has infinite
type $(a_{\theta},b_{\theta})$ with $a_{\theta}\not\equiv
b_{\theta}$ mod($\# \Y_K^*$). Suppose also that $p$ splits, and
that $(a_{\theta}-b_{\theta},p-1)=1$. Then $S=S'$.
\end{lem}
\begin{proof} Let $v$ be a prime of $K$ dividing $\f$. Let $v_0$ be
a prime of $K_0$ dividing $v$. Denote by $\Delta_{v_0}$ the
stabilizer of $v_0$ in $K_0$. We have then that $I_{v_0}\subset
\Delta_{v_0}\subset \Delta$ acts non-trivially in the Tate module
$T_pE$ via the Hecke character $\varphi_p$. Hence, $I_{v_0}$ acts
on $e_{\theta}(\otimes^w T_pE)$ via
$\psi_{\theta,p}|_{I_{v_0}}=\varphi_p^{a_{\theta}-b_{\theta}}$, as
$I_{v_0}$ lies in the kernel of the $p$-adic cyclotomic character
($v$ is prime to $p$). Since $p$ splits, we have that
$\#\Delta=(p-1)^2$, and since
 $(a_{\theta}-b_{\theta},p-1)=1$, $\varphi_p^{a_{\theta}-b_{\theta}}$
 acts non-trivial on $I_{v_0}$.
\end{proof}

\begin{lem} \label{4.6}
Suppose that $p$ splits in $K$ and suppose that $p-1\nmid
a_{\theta}+l+1$ or $p-1\nmid b_{\theta}+l+1$ or $p-1\nmid
a_{\theta}-b_{\theta}$. Then $\chi$, as $\Delta$-representation,
is not the cyclotomic representation.
\end{lem}
\begin{proof} The character $\chi$ is equal to
$(\overline{\psi}_{\theta}\kappa^l)^{-1}$ where $\kappa$ is the
cyclotomic character. Since $p$ is split in $K$, we have that
$p=\p\p^*$, with $\p\neq\p^*$. Let $\Delta_{\p}$ be the Galois
group $Gal(K(E[\p])/K)$; it is a subgroup of the decomposition
group since $\p$ is totally ramified in $\Delta_{\p}$. Observe
that $M_{\theta\Z_p}$ has multiplication isomorphic to
$\Y_K\otimes\Z_p$ and, as $p$ splits, it decomposes in two
idempotents. These idempotents decompose the Hecke character
$\psi_{\theta,p}=\psi_{\Omega_1}\oplus\psi_{\Omega_2}$, see
\cite{G} for more details. It is known that
$\overline{\psi}_{\Omega_1}|_{\Delta_{\p}}=\kappa^{b_{\theta}}$
(see for example \cite[\S 2.5]{G}), so we get that our character
is different from $\kappa$ as long as $\#\Delta_{\p}=p-1\nmid
b_{\theta}+l+1$, since $\kappa$ is a generator for the character
group of $\Delta_{\p}$.

Using the same kind of argument for $\p^*$ instead of $\p$ we
obtain a similar divisibility result but with $a_{\theta}$ instead
of $b_{\theta}$. Thus, we obtain the cyclotomic character only in
the case that $p-1\mid a_{\theta}+l+1$ and $p-1\mid
b_{\theta}+l+1$.

Similar argument for $\overline{\psi}_{\Omega_2}$, we obtain the
same simultaneous arithmetic conditions, i.e. $p-1|l+b_{\theta}+1$
and $p-1|a_{\theta}+l+1$ in order to obtain the cyclotomic
character. We refer to \cite[p.220,pp.223-234]{Ge} for more
details on the characters $\psi_{\Omega_1}$ and $\psi_{\Omega_2}$.
\end{proof}

Next theorem is the analog of \cite[Theorem 2.2.12]{K} in our
situation.

\begin{thm} Suppose that $p$ is an odd prime, prime to $N_{K/\Q}\mathfrak{f}$
and to $\#\Y_K^*$. Let the $\Delta$-representation $\chi$ on
$Hom_{\Y_p}(M_{\theta\Z_p}(w+l),\Y_p)$ be a good representation.
Then the map $(Soul)_p$ induces an isomorphism of $\Y_p$-modules
$$det_{\Y_p}((\overline{\Cc}_{\infty,\f_{\theta}}^{\chi}\otimes_{\Y_p}M_{\theta\Z_p}(w+l))
\otimes_{\Y_p[[\Gamma]]}^{\mathbb{L}}\Y_p)\cong
det_{\Y_p}(R\Gamma(\Y_S,M_{\theta\Z_p}(w+l+1)))^{-1}.$$
\end{thm}
\begin{proof}
In order to prove the theorem one can follow the same steps as in
the proof of \cite[Theorem 2.2.12]{K}, but with
$M_{\theta\Z_p}(w)$ instead of $T_pE$. The only results that need
some work are \cite[Proposition 2.2.13]{K} and \cite[Lemma
2.2.16]{K} (see \cite[\S 3.4]{Ba2} for a detailed proof). We will
show next the necessary steps to prove these two results.

Let us prove the analog of \cite[Proposition 2.2.13]{K}:
$$det_{\Y_p}(R\Gamma(\mathcal{G},H^0(K_{\infty}\otimes\Q_p,M_{\theta\Z_p}(w+l+1)'))^*)\cong\Y_p$$
and
$$det_{\Y_p}(R\Gamma(\mathcal{G},H^0(\Y_{\infty,S_p},M_{\theta\Z_p}(w+l+1)'))^*)\cong\Y_p,$$
where $M'$ denotes
$Hom_{\mathcal{O}_p}(M,\mathcal{O}_p\otimes\Q_p/\Z_p(1))$ and
$M^*=Hom_{\mathcal{O}_p}(M,\mathcal{O}_p\otimes\Q_p/\Z_p)$.

It follows from \cite[prop. 2.4.6]{G}, that the action of
$\mathcal{G}$ on $M_{\theta\Z_p}(w+l)$ is via the character
$$\overline{\psi}_{\theta,p}:\mathcal{G}\rightarrow(\Y_{\theta}\otimes\Z_p)^*,$$
multiplied by the $l$-th power of the $p$-adic cyclotomic
character.

Then it induces a surjection of $\Y_p$-modules
$\rho:\Y_p[[\Gamma]]\rightarrow M_{\theta\Z_p}(w+l)$ by the action
described above. Thus $ker(\rho)$ is an ideal of height 2 because
$\Gamma\cong\Z_p^2$. We know that $det_{R}$ is determined by the
ideals of height 1 for the ring $R$ (cf. \cite[2.1.4]{Ka1}). We
are going to show that this implies that
\begin{equation}\label{uiuiui}
det_{\Y_p}(M_{\theta\Z_p}(w+l)\otimes_{\Y_p[[\mathcal{G}]]}^{\mathbb{L}}\Y_p)\cong\Y_p.
\end{equation}
In fact, since $\Delta$ is finite and
$\mathcal{G}\cong\Gamma\times\Delta$, we have the isomorphism
$$M_{\theta\Z_p}(w+l)\otimes^{\mathbb{L}}_{\Y_p[[\mathcal{G}]]}\Y_p\cong
(M_{\theta\Z_p}(w+l))_{\Delta}\otimes^{\mathbb{L}}_{\Y_p[[\Gamma]]}\Y_p.$$
Since we know that $ker(\rho)$ has height 2, we have that
$det_{\Y_p[[\Gamma]]}((M_{\theta\Z_p}(w+l))_{\Delta})\cong\Y_p[[\Gamma]]$
and so
$det_{\Y_p}((M_{\theta\Z_p}(w+l))_{\Delta}\otimes^{\mathbb{L}}_{\Y_p[[\Gamma]]}\Y_p)\cong\Y_p$.
This shows (\ref{uiuiui}). We conclude by using \cite[Lemma
2.2.6]{K}.

Now, we show the analog of \cite[Lemma 2.2.16]{K}: the restriction
map induces isomorphisms
$$det_{\Y_p}(H^0(\Delta,R\Gamma(\Y_{0,S_p},M_{\theta\Z_p}(w+l+1)))\cong$$
$$det_{\Y_p}(H^0(\Delta,R\Gamma( \Y_{0,S},M_{\theta\Z_p}(w+l+1))))
\cong det_{\Y_p}(R\Gamma(\Y_S,M_{\theta\Z_p}(w+l+1))).$$

To show this consider the exact triangle
$$R\Gamma(\Y_{0,S_p},M_{\theta\Z_p}(w+l+1))\rightarrow R\Gamma(\Y_{0,S},M_{\theta\Z_p}(w+l+1))$$
$$\rightarrow \oplus_{v_0\in S\setminus S_p}R\Gamma_{k(v_0)}(\Y_v,M_{\theta\Z_p}(w+l))[1]$$
where $\Y_{v_0}$ is the local ring at $v_0$ and $S_p$ is the set
of places that divide $p$. Since $T_pE$ is unramified at the
places of $K_0$ in $S\setminus S_p$, the same is true for
$e_{\theta}(\otimes T_pE)(l+1)$. By the purity theorems in \'etale
cohomology we have that
$$R\Gamma_{k(v_0)}(\Y_{v_0},M_{\theta\Z_p}(w+l+1))\cong R\Gamma(k(v_0),M_{\theta\Z_p}(w+l))[-2].$$

It remains to prove only that $$H^0(\Delta,\oplus_{v_0\in
S\setminus S_p}R\Gamma(k(v_0),M_{\theta\Z_p}(w+l))=0.$$ To show
this result, observe that $$H^1(k(v_0),M_{\theta\Z_p}(w+l))\cong
M_{\theta\Z_p}(w+l)_{Gal(\overline{k(v_0)}/k(v_0))}$$ and
$H^0(k(v_0),M_{\theta\Z_p}(w+l))=0$ because $-w-2l\le -3$.

Now, let $v_0$ be a prime of $K_0$ dividing $v$ a prime of $K$
with $v|\f_{\theta}$ and let $\Delta_{v_0}$ be the stabilizer of
$v_0$. Since $I_{v_0}\subset\Delta_{v_0}$ acts non trivially on
the coinvariants
$M_{\theta\Z_p}(w+l)_{Gal(\overline{k(v_0)}/k(v_0))}$ because
$v_0\mid\f_{\theta}$, there are no fix elements.
\end{proof}

\section{The comparison between the map $r_p$ and $(Soul)_p$ in the
constructible $K$-elements}

Let's start recalling the result of Kings on the specialization of
the elliptic polylogarithm sheaf, which is an important key in his
proof of the Tamagawa number conjecture.

Let $E$ be an elliptic curve over a base scheme $T$, and denote by
$\overline{\pi}:E\rightarrow T$ the structural morphism, which is
proper and smooth. Consider $U=E\setminus e$, where $e$ is the
zero section of $E$. Consider the elliptic polylogarithm sheaf
$\Pol_{\Q_p}$ on $U$, which is a lisse pro-sheaf (i.e. a
projective limits of lisse sheaves) \cite[\S 3.2]{K}.

Let
$\mathcal{H}_{\Q_p}:=\underline{Hom}_{T}(R^1\overline{\pi}_*\Q_p,\Q_p)$.
Using $\Pol_{\Q_p}$ one defines the $p$-adic Eisenstein classes
associated to any integer $k$ and any $M$-torsion point $t$ in $E$
different from $e$ as some elements in
$H^1(T,Sym^k\mathcal{H}_{\Q_p})$. The definition is extended by
linearity to any divisor supported on $M$-torsion points
(\cite[Def. 3.5.9]{K}). The main part of the result of Kings is
the explicit computation of these Eisenstein classes.


We are going to explain this result. Consider $H_n:=ker[p^n]$ as a
scheme over $T$. Let us consider the map multiplication by $p^n$,
$p_n:E_n\rightarrow E$, where $E_n$ is the elliptic curve $E$ over
$T$ considered as a $H_n$-torsor over $E$. Consider the
characteristic group $I[H_n]:=ker(p_{n,*}\Z\rightarrow\Z)$, which
is the characteristic group of a torus $T_{H_n}$. In this
situation we have the connecting map $\delta$ from the Kummer
exact sequence:
\begin{equation}\label{above}
\delta:H^0(H_n,T_{H_n})\rightarrow H^1(H_n,T_{H_n}[p^r]).
\end{equation}
Using this connecting morphism, we can express the Eisenstein
classes explicitly.
\begin{thm}[Kings, theorem 4.2.9 in \cite{K}]\label{poleli} Let $p$ be a prime number, and
let $E$ be an elliptic curve over a base scheme $T$ where $p$ is
invertible.

Let $\beta$ be any divisor in $E$ of the form
$$\beta:=\sum_{t\in E[M](T)\setminus e}n_t(t),$$
$n_t$ an integer and consider $[\mathfrak{a}]:E\rightarrow E$ any
isogeny with degree prime to $Mp$.

Then, for any $m>0$, the $p$-adic Eisenstein class
$$N\mathfrak{a}(\mathfrak{a}^{\otimes m}N\mathfrak{a}-1)(\beta^*Pol_{\Q_p})^m\in H^1(T, Sym^m\mathcal{H}_{\Q_p}(1))$$ is given by
$$\pm\frac{1}{m!}(\delta \sum_{t\in E[M](T)\setminus e}n_t\sum_{[p^n]t_n=t}\theta_{\mathfrak{a}}(-t_n)\widetilde{t_n}^{\otimes m})_n$$
where $\widetilde{t_n}$ is the projection of $t_n$ to $E[p^n]$ and
$\delta$ is the Sym-extension of the boundary map
$H^0(H_n,T_{H_n})\rightarrow H^1(H_n,T_{H_n}[p^r])$ where
$H_n:=ker[p^n]$ is considered as a scheme over $T$ and $T_{H_n}$
is the torus with character group $I[H_n]:=ker(p_{n,
*}\Z\rightarrow\Z)$.
\end{thm}

The following result relates the image of
$\mathcal{E}^{m}_{\Mm}(\beta)$ by the Soul\'{e} regulator map with
the polylogarithmic sheaf, where $\mathcal{E}^m_{\Mm}$ is the
Eisenstein symbol \cite[\S1.2.1]{K}.

\begin{thm}\label{polelidos} Under the same hypothesis of Theorem \ref{poleli}, let $\beta$ be as in the previous theorem. Then
$$r_p(\mathcal{E}^{m}_{\Mm}(\beta))=-M^{2m}(\beta^*Pol_{\Q_p})^{m}$$ in $H^1(T,Sym^{m}\mathcal{H}_{\Q_p}(1))$.
\end{thm}
\begin{proof} The same proof of \cite[Theorem 1.2.5]{K} with $m$
instead of $2k+1$ works. See also \cite[proof Theorem 3.5.2]{Ba2}.
\end{proof}

We are going to apply these results to the divisor
 $\beta_{\theta}=N_{K(\f)/K}((t))$, where $t:=\Omega f^{-1}$ is a $\f_{\theta}$-torsion point.
Take $M=N_{K/\Q}\f_{\theta}$, $m=w+2l$, $T=\Y_S$ and
$\mathcal{H}_{\Q_p}=T_pE\otimes\Q_p$, using the notations of the
theorem \ref{poleli}. Let $\mathfrak{a}\subset\Y_K$ be an ideal
prime to $6p\mathfrak{f}$, and consider the isogeny given by
$\varphi(\mathfrak{a})$. Let $\theta_{\mathfrak{a}}$ be the
classical theta function.

To simplify the notation, define for any $\widetilde{t_r}\in
E[p^r]$
$$\gamma(\widetilde{t_r})^m:=<\widetilde{t_r},\sqrt{d_K}\widetilde{t_r}>^{\otimes m}$$
where $<,>$ denotes the Weil pairing. Our objective is the
computation of
$$\mathcal{K}_{\Mm}\circ\mathcal{E}^{w+2l}_{\Mm}(\beta_{\theta}).$$
Remember that we are under the restriction $a_{\theta}\not\equiv
b_{\theta}\ mod(\# \Y_K^*)$.

We consider the following commutative diagram \cite[(2.8)]{De}
\begin{small}
\begin{center}
\begin{tabular}{lcr}
$H^{2l+w+1}_{\Mm}(Sym^{2l+w}h^1E,\Q(w+2l+1))$&$\overset{((\Delta_{CM})^l\times id)^{*}}{\longrightarrow}$&$H^{2l+w+1}_{\Mm}(E^{l+w},\Q(2l+w+1))$\\
$\mbox{}\ \mbox{} \ \mbox{} \Kk_{\Mm}\downarrow$&&$\downarrow pr_{*}\mbox{}\ \mbox{} \ \mbox{} $\\
$H^{w+1}_{\Mm}(M_{\theta\Q},\Q(w+l+1))$&$\overset{e_{\theta}}{\longleftarrow}$&$H^{w+1}_{\Mm}(h^1(E)^{\otimes w},\Q(l+w+1)),$\\
\end{tabular}
\end{center}
\end{small}
where $pr$ is the projection in the last $w$ components and
$\Delta_{CM}:E\rightarrow E\times E$ is given by
$e\mapsto(e,\sqrt{d_K}e)$. We obtain a map in Galois cohomology
given by
$$H^1(\Y_S,Sym^{2l+w}\mathcal{H}_{\Q_p}(1))\rightarrow $$
$$H^1(\Y_S,(e_{\theta}\mathcal Sym^w{\mathcal{H}}_{\Q_p})(l+1))=H^1(\Y_S,M_{\theta\Q_p}(w+l+1))$$
such that
$$\mathcal{K}_{\Mm}(\varphi_p(\mathfrak{a})^{\otimes 2l+w}
Sym^{2l+w}\mathcal{H}_{\Q_p}(1))=e_{\theta}(\otimes^w\varphi_p(\mathfrak{a}))N\mathfrak{a}^lSym^w\mathcal{H}_{\Q_p}(l+1).$$
\begin{thm} \label{aii:teo}
Let $p$ be a prime number such that $p\nmid
6N_{K/\Q}(\mathfrak{f})$. Let $\theta$ be an idempotent with
infinity type $(a_{\theta},b_{\theta})$ which
$a_{\theta}\not\equiv b_{\theta}\ mod\ (\#\Y_K^*)$. For a
$p^rN_{K/\Q}(\mathfrak{f}_{\theta})$-torsion point $t_r$, denote
by $\widetilde{t_r}$ its projection to $E[p^r]$. Then, if
$t=\Omega f^{-1}$, we have the following equality
$$N\mathfrak{a}\left(\psi_{\theta,p}(\mathfrak{a})N\mathfrak{a}^{l+1}-1\right)
r_p(\xi_{\theta,l})=$$
\begin{small}
$$\frac{(-1)^{l}L_p(\overline{\psi}_{\theta},-l)^{-1}
N_{T_{\theta}/\Q}\mathfrak{f}_{\theta}^{3l+2w}\Phi(\f_{\theta})
}{2^{l-1} \psi_{\theta}(\rho_{\theta})\Phi(\f)} \cdot ( \left(
\delta
N_{K(\mathfrak{f})/K}\sum_{p^rt_r=t}\theta_{\mathfrak{a}}(-t_r)\otimes
e_{\theta} (\otimes^w\widetilde{t_r})\otimes
\gamma(\widetilde{t_r})^l\right)_r$$
\end{small}

\end{thm}
\begin{proof} Using theorems \ref{yok:teo}, \ref{poleli} and \ref{polelidos}, we have that
$$r_p(\xi_{\theta,l})=\frac{(-1)^{l-1}(2l+w)!L_p(\overline{\psi}_{\theta},-l)^{-1}\Phi(\f_{\theta})
}{2^{l-1}N_{T_{\theta}/\Q}\mathfrak{f}_{\theta}^l
\psi_{\theta}(\rho_{\theta})\Phi(\f)}\mathcal{K}_{\Mm}(\mathcal{E}_{\Mm}^{2l+w}(\beta))$$
$$=\frac{(-1)^{l}(2l+w)!L_p(\overline{\psi}_{\theta},-l)^{-1}
N_{T_{\theta}/\Q}\mathfrak{f}_{\theta}^{3l+2w}\Phi(\f_{\theta})
}{2^{l-1}
\psi_{\theta}(\rho_{\theta})\Phi(\f)}\mathcal{K}_{\Mm}(\beta^*Pol_{\Q_p})^{w+2l}.$$
We have that
$$\mathcal{K}_{\Mm}(\widetilde{t_r}^{\otimes
2l+w})=e_{\theta}(\otimes^w\widetilde{t_r})\otimes
\gamma(\widetilde{t_r})^l.$$ Finally, applying Kings' theorem
\ref{poleli}, we obtain the desired identity.
\end{proof}
We want to rewrite the previous formula in terms of the norm map
of the extension $K(\f_{\theta})K(E[p^n])/K$. We will work with
$\f$ instead of $\f_{\theta}$ since then we can use that
$K(E[\p^n\mathfrak{f}])=K(\mathfrak{p^nf})$, the ray class field,
because $\f$ is the conductor of $E$ and divides the ideal
$\f\p^n$ (\cite[II, Prop.1.6]{dS}).

Fix a prime $\mathfrak{p}$ of $K$ where $E$ has good reduction,
and take $\pi=\varphi(\mathfrak{p})$. Denote by
$$H^{\mathfrak{p}}_{r,t}:=\{ t_r\in E[\mathfrak{p}^r\mathfrak{f}]|\pi^rt_r=t\}.$$
Let $\sigma_{\mathfrak{p}}$ be the Frobenius at $\mathfrak{p}$ in
$Gal(K(\f)/K)$, we write $t_r=(\widetilde{t_r},\pi^{-r}t)\in
E[\mathfrak{p}^r\mathfrak{f}]=E[\mathfrak{p}^r]\oplus
E[\mathfrak{f}]$, where $\pi^{-r}t$ means
$t^{\sigma_{\mathfrak{p}}^{-r}}$. \\ Consider the filtration of
$H^{\mathfrak{p}}_{r,s}$ defined by
$$F^i_{r,t}:=\{t_r\in H^{\mathfrak{p}}_{r,s}|\pi^{r-i}\widetilde{t_r}=0\}.$$

\begin{thm} Let $\mathfrak{p}$ be as above and $t_r=(\widetilde{t_r},\pi^{-r}t)\in F^0_{r,s}\setminus F^1_{s,t}$.
Suppose that $\Y_K^*\rightarrow(\Y_K/\f_{\theta})^*$ is injective.
Denote the Euler factor of the Hecke character
$\overline{\psi}_{\theta}$ at $\p$ evaluated at $-l$ by
$L_{\mathfrak{p}}(\overline{\psi_{\theta}},-l)$. Then
$$L_{\mathfrak{p}}(\overline{\psi}_{\theta},-l)^{-1}\left(N_{K(\mathfrak{f})/K}
\sum_{s_r\in
H^{\mathfrak{p}}_{r,t}}\theta_{\mathfrak{a}}(-s_r)\otimes
e_{\theta}(\otimes^w
\widetilde{s_r})\otimes\gamma(\widetilde{s_r})^l\right)_r=$$
$$(N_{K(\mathfrak{p}^r\mathfrak{f})/K}\left(\theta_{\mathfrak{a}}(-t_r)\otimes e_{\theta}(\otimes^w\widetilde{t_r})\otimes\gamma(\widetilde{t_r})^l\right))_r$$
in $H^1(\Y_S,e_{\theta}(T_{\mathfrak{p}}E(1))(l)\otimes\Q_p)$ for
all $\mathfrak{a}$ relatively prime to $\mathfrak{pf}$.
\end{thm}
\begin{proof} The identification $Hom_{\Y_p}(T_pE,\Y_p)\cong T_pE(-1)$ is via the
conjugate linear $\Y_p$-action on the right side. Hence
$\overline{\varphi(\mathfrak{p})}t_r=t_{r-1}$. We have the
equality
$$(\overline{\psi_{\theta}(\mathfrak{p})}/N\mathfrak{p}^{-l})^i
N_{K(\mathfrak{p}^r\mathfrak{f})/K(\mathfrak{p}^{r-i}\mathfrak{f})}(\theta_{\mathfrak{a}}(-t_r)\otimes
e_{\theta}(\otimes^w\widetilde{t_r})\otimes\gamma(\widetilde{t_r})^l)=$$
$$N_{K(\mathfrak{p}^r\mathfrak{f})/K(\mathfrak{p}^{r-i}\mathfrak{f})}(\theta_{\mathfrak{a}}(-t_r)\otimes
e_{\theta}(\otimes^w\overline{\varphi(\mathfrak{p})}^i
\widetilde{t_r})\otimes\gamma(\overline{\varphi(\mathfrak{p})}^i\widetilde{t_r})^l)=$$
$$(N_{K(\mathfrak{p}^r\mathfrak{f})/K(\mathfrak{p}^{r-i}\mathfrak{f})}(\theta_{\mathfrak{a}}(-t_r)))\otimes
e_{\theta}(\otimes^w\widetilde{t_{r-i}})\otimes\gamma(\widetilde{t_{r-1}})^l)=$$
$$\theta_{\mathfrak{a}}(-(\widetilde{t_{r-i}},\pi^{i-r}t))\otimes
e_{\theta}(\otimes^w
\widetilde{t_{r-i}})\otimes\gamma(\widetilde{t_{r-i}})^l,$$ where
the last equality uses the distribution relation for
$\theta_{\mathfrak{a}}$ (\cite[II 2.5]{dS}).

The Galois group of
$K(\mathfrak{p}^{r-i}\mathfrak{f})/K(\mathfrak{f})$ acts simply
transitively on $F^i_{r,t}\setminus F^{i+1}_{r,t}$. We get that
$$(\overline{\psi_{\theta}(\mathfrak{p})}/N\mathfrak{p}^{-l})^i
N_{K(\mathfrak{p}^r
\mathfrak{f})/K(\mathfrak{f})}(\theta_{\mathfrak{a}}(-t_r)\otimes
e_{\theta}
(\otimes^w\widetilde{t_r})\otimes\gamma(\widetilde{t_r})^l)=$$
$$\sum_{t_{r-i}\in F^i_{r,t}\setminus
F^{i+1}_{r,t}}\theta_{\mathfrak{a}}(-(\widetilde{t_{r-i}},\pi^{i-r}t))\otimes
e_{\theta}(\otimes^w
\widetilde{t_{r-i}})\otimes\gamma(\widetilde{t_{r-i}})^l.$$

We know by \cite[Prop. II.2.4.ii)]{dS} that we have the equality
$\theta_{\mathfrak{a}}(-(\widetilde{t_{r-i}},\pi^{i-r}t))=\theta_{\mathfrak{a}}(-(\widetilde{t_{r-i}},\pi^{-r}t))^{\sigma^i_{\mathfrak{p}}}$
with $\sigma_{\mathfrak{p}}$ is the Frobenius at $\mathfrak{p}$ in
the Galois group of $K(\mathfrak{f})/K$. This and the fact that
$N_{K(\mathfrak{f})/K}$ is the sum over all Galois translates,
which act trivially on $\widetilde{t_{r-i}}$, gives that
$$(\overline{\psi_{\theta}(\mathfrak{p})}/N\mathfrak{p}^{-l})^i
N_{K(\mathfrak{p}^r\mathfrak{f})/K}(\theta_{\mathfrak{a}}(-t_r)\otimes
e_{\theta}(\otimes^w\widetilde{t_r})\otimes\gamma(\widetilde{t_r})^l)=$$
$$N_{K(\mathfrak{f})/K}\left( \sum_{t_{r-i}\in F^i_{r,t}\setminus
F^{i+1}_{r,t}}
\theta_{\mathfrak{a}}(-(\widetilde{t_{r-i}},\pi^{-r}t))\otimes
e_{\theta}(\otimes^w
\widetilde{t_{r-i}})\otimes\gamma(\widetilde{t_{r-i}})^l\right),$$
Adding these equalities with respect to $i$ and increasing $r$ if
necessary we get the result.
\end{proof}
\begin{lem} \label{cincsis}\label{5.5} Suppose that $\theta$ has infinity type $(w,0)$ or $(0,w)$ and $(\#\Y_K^*,w)=1$. Then
$$\Y_K^*\rightarrow (\Y_K/\f_{\theta})^*$$ is injective.
\end{lem}
\begin{proof} Let $u$ be and element in $\Y_K^*$, $u\neq1$ and consider the id\`ele defined by
$x_{\infty}=1$ and $x_{\p}=u$ at all finite places $\p$ of $K$.
Then as complex Hecke character
$\varphi^w(x)=\varphi^w(u^{-1}x)=u^w\neq1$ if $(w,\#\Y_K^*)=1$.
So, by definition of the conductor of $\psi_{\theta}$, we obtain
that $u\not\equiv 1 (mod\ \f_{\theta})$, hence the result for the
type $(w,0)$.
 For the type $(0,w)$ the
proof is similar but with $\overline{\varphi}$ instead of
$\varphi$.
\end{proof}

\begin{cor}\label{nor:teo} Suppose $p\nmid 6N_{K/\Q}(\f)$, $\Y_K^*\rightarrow(\Y_K/\f_{\theta})^*$
is injective and the hypothesis of theorem \ref{aii:teo}. Then
$$N\mathfrak{a}(\psi_{\theta,p}(\mathfrak{a})N\mathfrak{a}^{l+1}-1)r_p(\xi_{\theta,l})=$$
$$\pm\frac{N_{K/\Q}(\mathfrak{f}_{\theta})^{3l+2w}\Phi(\f_{\theta})}{2^{l-1}\psi_{\theta}(\rho_{\theta})\Phi(\f)}
\delta\left(
N_{K(E[p^r])K(\mathfrak{f})/K}\theta_{\mathfrak{a}}(-t_r)\otimes
e_{\theta}(\otimes^w\widetilde{t_r})\otimes\gamma(\widetilde{t_r})^l\right)_r=$$
$$\pm\frac{N_{K\Q}(\mathfrak{f}_{\theta})^{3l+2w}}{2^{l-1}\psi_{\theta}(\rho_{\theta})}
\cdot
 \delta\left(
N_{K(E[p^r])K(\mathfrak{f}_{\theta})/K}\theta_{\mathfrak{a}}(-t_r)\otimes
e_{\theta}(\otimes^w\widetilde{t_r})\otimes\gamma(\widetilde{t_r})^l\right)_r$$
where $t_r$ is a primitive $p^r\mathfrak{f}_{\theta}$-torsion
point with $p^rt_r=t$ and $\mathfrak{a}$ is relative prime to $p
\f$.
\end{cor}

\begin{proof} If $p$ is inert or prime the first equality is deduced from the previous theorem.
If $p$ split, it decomposes in a $\p$ part and a $\p^*$ part.
Putting together the previous result with $\p$ and with $\p^*$, we
have the first equality.

To show the second equality, consider
$$N_{K(E[p^r])K(\f)/K(E[p^r])K(\f_{\theta})}\theta_{\mathfrak{a}}(-t_r)=$$
$$\prod_{\sigma\in Gal(K(\f)K(E[p^r])K(p^r)/K(\f_{\theta})K(E[p^r])K(p^r))}\theta_{\mathfrak{a}}(-t_r)^{\sigma}$$
because $K(\f)$ is disjoint with $K(p^r)$ over $K$ since $K=K(1)$,
and $K(\f)=K(E[\f])$ is disjoint with $K(E[p^r])$ over $K$.
Moreover, since $\theta_{\mathfrak{a}}(-t_r)\in K(\f)K(p^r)=K(\f
p^r)$ and $(\f,p)=1$, we have that the norm is equal to
$$\prod_{\tau\in Gal(K(\f p^r)/K(\f_{\theta}p^r))}\theta_{\mathfrak{a}}(-t_r)^{\tau}.$$
But $\theta_{\mathfrak{a}}(-t_r)\in K(\f_{\theta}p^r)$ because
$-t_r$ is a point of $\f_{\theta}p^r$-torsion. To obtain the
second equality we need to show that
$\frac{\Phi(\f_{\theta})}{\Phi(\f)}\#(Gal(K(\f)/K(\f_{\theta})))$
is one. We have
$\frac{\Phi(\f_{\theta})}{\Phi(\f)}\#(Gal(K(\f)/K(\f_{\theta})))=\frac{\Phi(\f_{\theta})}{\Phi(\f)}\frac{[K(\f):K(1)]}{[K(f_{\theta}):K(1)]},$
we are in Galois extensions. Observe that
$\frac{\Phi(\f_{\theta})}{\Phi(\f)}\frac{[K(\f):K(1)]}{[K(f_{\theta}):K(1)]}
=\frac{\Phi(\f_{\theta})}{\Phi(\f)}\frac{\Phi(\f)\frac{w_{\f}}{w_K}}{\Phi(\f_{\theta})\frac{w_{\f_{\theta}}}{w_K}}$
where $w_K$ are the roots of unity of $K$ and $w_{\mathfrak{m}}$
are the roots of unit of $\mathcal{O}_K^*$ congruent to 1 modulo
$\mathfrak{m}$, by class field theory (see \cite[p.36]{dS}). And
is equal to one because $\mathcal{O}^*\rightarrow
(\Y_K/\mathfrak{m})^*$ is injective for $\mathfrak{m}=\f_{\theta}$
by hypothesis, and for $\mathfrak{m}=\f$ (Lemma \ref{cincsis} or
\cite[Corollary 5.9]{Ru4}).
\end{proof}

Now we want to show that the elements
$$(N_{K(E[p^r])K(\mathfrak{f}_{\theta})/K} \theta_{\mathfrak{a}}(t_r)\otimes
e_{\theta}(\otimes^w\widetilde{t_r})\otimes\gamma(\widetilde{t_r})^l)_r$$
generate $(\overline{\Cc}_{\infty}^{\chi}\otimes
M_{\theta\Z_p}(w+l))_{\Gamma}$, where $\mathfrak{a}$ is prime to
$6p\mathfrak{f}$ and $\chi$ is the representation of $\Delta$ on
$Hom_{\Y_p}(M_{\theta\Z_p}(w+l),\Y_p)$, that we suppose a good
representation.

We suppose from now on that the natural map
$\Y_K^*\rightarrow(\Y_K/\f_{\theta})^*$ is injective, assumption
also needed to define our elliptic units.

\begin{prop}\label{gen:teo} Consider $p\nmid 6N_{K/\Q}(\mathfrak{f})$ and $\mathfrak{a}$
an ideal in $\Y_p$, which is prime to $6p\mathfrak{f}$ and such
that $\psi_{\theta,p}(\mathfrak{a})N\mathfrak{a}^{l+1}\not\equiv
1(mod\ p)$. Then the $\Y_p[[\Gamma]]$-module
$$\overline{\Cc}_{\infty,\f_{\theta}}^{\chi}\otimes_{\Y_p}M_{\theta\Z_p}(w+l)$$ is generated by $(\theta_{\mathfrak{a}}(t_r)\otimes
e_{\theta}(\otimes^w\widetilde{t_r})\otimes\gamma(\widetilde{t_r})^l)_r$,
where $t_r$ is a primitive $p^r\mathfrak{f}_{\theta}$-division
point.
\end{prop}

\begin{remark} The existence of
 an ideal $\mathfrak{a}$ satisfying the conditions of the proposition \ref{gen:teo} is
 equivalent to the condition that the $\Delta$-representation $\chi$
 is not the cyclotomic representation.
\end{remark}

\begin{proof}Observe first that
$e_{\theta}(\widetilde{t_r})\otimes\gamma(\widetilde{t_r})$
generates $M_{\theta\Z_p}(w+l)$, because $M_{\theta\Z_p}(w)$ is
one dimensional and concerning how it generates $\Z_p(l)$ use the
same proof did in \cite[p.623]{K}.

Remember that we have an inclusion of
$\overline{\Cc}^{\chi}_{\infty,\f_{\theta}}$ in $\Uu_{\infty}$ the
local units Iwasawa module, which is torsion free
\cite[Prop.11.4]{Ru4}, thus
$\overline{\Cc}^{\chi}_{\infty,\f_{\theta}}$ is a torsion free
$\Y_p[[\Gamma]]$-module. Is enough to show that is one
dimensional. Let $\mathfrak{b}$ be another ideal prime to
$6p\mathfrak{f}$. Take
$\sigma_{\mathfrak{a}}=[\mathfrak{a},K_n/K]$ and
$\sigma_{\mathfrak{b}}=[\mathfrak{b},K_n/K]$. Then, by the
properties of the theta function, we have that
$$(\sigma_{\mathfrak{a}}-\psi_{\theta,p}(\mathfrak{a})N\mathfrak{a}^{l+1})
(\theta_{\mathfrak{b}}(t_n)\otimes
e_{\theta}(\otimes^w\widetilde{t_n})\otimes
\gamma(\widetilde{t_n})^l)=$$
$$\psi_{\theta,p}(\mathfrak{a})N\mathfrak{a}^l(\theta_{\mathfrak{b}}(t_n)^{\sigma_{\mathfrak{a}}-N\mathfrak{a}}\otimes
e_{\theta}(\otimes^w\widetilde{t_n})\otimes\gamma(\widetilde{t_n})^l)=$$
$$\psi_{\theta,p}(\mathfrak{a})N\mathfrak{a}^l(\theta_{\mathfrak{a}}(t_n)^{\sigma_{\mathfrak{b}}-N\mathfrak{b}}\otimes
e_{\theta}(\otimes^w\widetilde{t_n})\otimes\gamma(\widetilde{t_n})^l).$$
 Then, it is enough show that
$(\sigma_{\mathfrak{a}}-\psi_{\theta,p}(\mathfrak{a})N\mathfrak{a}^{l+1})$
is invertible in $\Y_p[[\Gamma]]$. But the element
$\sigma_{\mathfrak{a}}$ corresponds to 1 on $\Y_p/p$ and thus
$\sigma_{\mathfrak{a}}-\psi_{\theta}(\mathfrak{a})N\mathfrak{a}^{l+1}$
is invertible in $\Y_p[[\Gamma]]$ because
$1\not\equiv\psi_{\theta}(\mathfrak{a})N\mathfrak{a}^{l+1}\ mod\
p$.
\end{proof}

\begin{cor}\label{cor5.9} Assume that $p\nmid 6N_{K/\Q}\f$. Then the image of $\Rr_{\theta}$ by $r_p$ in
the cohomology group $H^1(\Y_S,M_{\theta\Z_p}(w+l+1))\otimes\Q_p$
coincides with
$$(Soul)_p((\overline{\Cc}^{\chi}_{\infty,\f_{\theta}}\otimes M_{\theta\Z_p}(w+l))_{\Gamma}).$$
\end{cor}
\begin{proof} As
$$N\mathfrak{f}_{\theta}^{3l+2w}/2^{l-1}\psi_{\theta}(\rho_{\theta})$$
is prime to $p$, it follows from the definition of $(Soul)_p$ and
Corollary \ref{nor:teo}.
\end{proof}

\begin{lem}
The canonical map
 $$(\overline{\Cc}_{\infty,\f_{\theta}}\otimes M_{\theta\Z_p}(w+l))
 \otimes^{\mathbb{L}}_{\Y_p[[\mathcal{G}]]}\Y_p\rightarrow
 (\overline{\Cc}_{\infty,\f_{\theta}}\otimes M_{\theta\Z_p}(w+l))_{\mathcal{G}}\cong
 (\overline{\Cc}_{\infty,\f_{\theta}}^{\chi}\otimes M_{\theta\Z_p}(w+l))_{\Gamma}$$
is an isomorphism and moreover
$(\overline{\Cc}_{\infty,\f_{\theta}}^{\chi}\otimes
M_{\theta\Z_p}(w+l))_{\Gamma}\cong\Y_p$.

\end{lem}
\begin{proof} We observe that the proof of proposition \ref{gen:teo} shows that
$\overline{\Cc}^{\chi}_{\infty,\f_{\theta}}\cong\Y_p[[\Gamma]]$ is
a free $\Y_p[[\Gamma]]$-module of rank 1. This implies, as in
\cite[lemma 5.2.3]{K}, that
$(\overline{\Cc}^{\chi}_{\infty,\f_{\theta}}\otimes
M_{\theta\Z_p}(w+l))_{\Gamma}\cong\Y_p$. The claim follows since
the previous module is induced and hence the higher Tor-terms
vanish.
\end{proof}

\begin{cor}\label{end:teo} The map
$$\Rr_{\theta}\otimes\Z_p\rightarrow R\Gamma(\Y_S,M_{\theta\Z_p}(w+l+1)\otimes\Q_p)[1]$$
induced by $r_p$, gives an isomorphism
$$det_{\Y_p}\Rr_{\theta}\cong det_{\Y_p}R\Gamma(\Y_S,M_{\theta\Z_p}(w+l+1))^{-1}$$
\end{cor}

Before stating the next theorem, let us recall all the hypothesis
we used during the paper and that we will need:\\ {\bf
($\diamond\diamond\diamond$)}
 Let $p$ be a fix prime such that $p\nmid 6N_{K/\Q}\f$ (hence, in particular
 $p\nmid (\#\Y_K^*)$, and $p\nmid D_K$). Consider
$l$ a non-negative integer. Let $(a_{\theta},b_{\theta})$ be the
infinite type of $\psi_{\theta}$, with $a_{\theta},b_{\theta}$
non-negative integers with $w=a_{\theta}+b_{\theta}\ge1$ such that
$a_{\theta}\not\equiv b_{\theta}\ mod(\#\Y_K^*)$ and $-w-2l\le-3$.
Assume that $\Y_K^*\rightarrow(\Y_K/\f_{\theta})^*$ is injective.
Suppose moreover that the representation $\chi$ of $\Delta$ in
$Hom_{\Y_p}(H^w_{\acute{e}t}(M_{\theta}\times_K\overline{K},\Z_p(w+l)),\Y_p)$
is a good representation (see the definition in \ref{goo:teo})
which is not equal as $\Delta$-representation to the cyclotomic
representation.

\begin{thm}\label{kht:teo} Under the hypothesis {\rm {\bf ($\diamond\diamond\diamond$)}} above,
 there is an $\Y_K$-submodule $\Rr_{\theta}\subset
H^{w+1}_{\Mm}(M_{\theta\Q},\Q(w+l+1))$ of rank 1 such that:
\begin{enumerate}
\item
$det_{\Y_{K}[1/D_K]}(r_{\Dd}(\Rr_{\theta}\otimes_{\Y_K}\Y_K[1/D_K]))\cong$
$$L_S^*(\overline{\psi}_{\theta},-l)
det_{\Y_K[1/D_K]}(H^w_B(M_{\theta\C},\Z(w+l))\otimes_{\Y_K}\Y_K[1/D_K])$$
in
$det_{\Y_K[1/D_K]\otimes\R}(H_B^w(M_{\theta\C},\Z(w+l))\otimes_{\Y_K}\Y_K[1/D_K]\otimes\R)$.
\item The map $r_p$ induces an isomorphism
$$det_{\Y_{K}\otimes\Z_p}(\Rr_{\theta})\cong det_{\Y_{K}\otimes\Z_p}(R\Gamma(\Y_{K}[1/S],M_{\theta\Z_p}(w+l+1))^{-1}.$$
\end{enumerate}
Here
$$L_S^*(\overline{\psi}_{\theta},-l)=\lim_{s\rightarrow-l}\frac{L_S(\overline{\psi}_{\theta},s)}{s+l},$$
and $S$ is the set of primes of $K$ dividing $p$ and the ones
dividing $\f_{\theta}$.

Moreover, if $r_p$ is injective on $\Rr_{\theta}$, the second part
can be written as
$$det_{\Y_{K}\otimes\Z_p}(H^1(\Y_{K}[1/S],M_{\theta\Z_p}(w+l+1))/r_p(\Rr_{\theta}))\cong $$
$$det_{\Y_{K}\otimes\Z_p} H^2(\Y_{K}[1/S],M_{\theta\Z_p}(w+l+1)).$$

\end{thm}

\begin{proof} It is a direct consequence of the theorem \ref{yok:teo} and the above
corollary \ref{end:teo}.
\end{proof}

After taking the norm $N_{K/\Q}$, we obtain the following result.

\begin{thm}\label{qht:teo} Under the assumption {\rm {\bf ($\diamond\diamond\diamond$)}} above,
 there is a $\Z$-submodule $\Rr_{\theta}$ in
$H^{w+1}_{\Mm}(M_{\theta\Q},\Q(w+l+1))$ of rank 2 such that:
\begin{enumerate}
\item The map $r_{\Dd}\otimes\R$ is an isomorphism restricted to
$\Rr_{\theta}\otimes\R$. \item
$dim_{\Q}(H^w_B(M_{\theta\C},\Q(w+l)))=ord_{s=-l}L_S(M_{\theta\Q},s)=2$.
\item We have the equality
$$r_{\Dd}(det_{\Z[1/D_K]}(\Rr_{\theta}\otimes_{\Y_K}\Y_K[1/D_K]))=$$
$$L_S^*(M_{\theta\Q},-l)det_{\Z[1/D_K]}(H^w_B(M_{\theta\C},\Z(w+l))\otimes_{\Y_K}\Y_K[1/D_K])$$
where
$$L^*_S(M_{\theta\Q},-l)=\lim_{s\rightarrow-l}\frac{L_S(M_{\theta\Q},s)}{(s+l)^2}$$
and $S$ is the set of places of $K$ that divides $p$ and the
places dividing the conductor $\f_{\theta}$. \item We have that
$$det_{\Z_p}(\Rr_{\theta}\otimes \Z_p)=det_{\Z_p}(R\Gamma(\Y_{K}[1/S],M_{\theta\Z_p}(w+l+1)))^{-1}.$$
If $r_p$ is injective on $\Rr_{\theta}$, then
$r_p(det_{\Z_p}(\Rr_{\theta}\otimes\Z_p))$ is a basis of the
$\Z_p$-lattice
$$det_{\Z_p}(R\Gamma(\Y_{K}[1/S],M_{\theta\Z_p}(w+l+1)))^{-1}$$
$$\subset det_{\Q_p}(R\Gamma(\Y_{K}[1/S],M_{\theta\Z_p}(w+l+1)\otimes\Q)[-1]).$$
\end{enumerate}

\end{thm}

\begin{remark}\label{5.14} Theorems \ref{kht:teo} and \ref{qht:teo} imply the weak
$p$-part of the Tamagawa number conjecture for Hecke characters
\cite{Ka1} for $K$ or $\Q$ coefficients respectively, up to the
finiteness of $H^2_p:=H^2(\Y[1/S],M_{\theta\Z_p}(w+l+1))$ and the
bijectively of the Soul\'e regulator map $r_p$. Concerning these
requirements, we have the following.
\begin{enumerate} \item If $p$ is a regular prime for the field $K(E[p])$, then
$H^2_p$ is finite \cite{Ba3}. Moreover without any assumption, one
obtains that for almost all $l$ this Galois cohomology group is
finite \cite[Theorem 12.4]{Ka6} or \cite{Ba4}.

\item About the bijectively of the Soul\'e regulator map observe
if $H^2_p$ is finite, similar arguments as in \cite[\S 5.2.2]{K}
implies the injectivity for $(Soul)_p$ and therefore $r_p$ is
injective on $\Rr_{\theta}\otimes\Q_p$ by corollary \ref{cor5.9}.
Therefore $r_p$ restricted to $\Rr_{\theta}\otimes\Q_p$ is an
isomorphism \cite[cor. 1]{Ja}.
\end{enumerate}
Therefore for regular primes $p$, we obtain in full generality the
weak $p$-part of the Tamagawa number conjecture for Hecke
characters of imaginary quadratic fields.
\end{remark}

\section{The remaining Tate twists}

\begin{subsection}{The remaining non-critical twists}

The value of the $L$-function at zero for $M_{\theta}(w+l+1)$ with
$-w-2l<-2$ is related with the first non-zero coefficient of the
Taylor development at $-l$ of the $L$-function associated to
$\overline{\psi}_{\theta}$ by the use of the functional equation
of $L$-functions. The non-critical values associated to the Hecke
character $\overline{\psi}_{\theta}$ (we restrict to the situation
$a_{\theta}\not\equiv b_{\theta}(mod|\Y_K^*|)$)
are the integers $l$ such that $-l\leq min(a_{\theta},b_{\theta})$
where $a_{\theta},b_{\theta}$ are associated to the Hecke
character $\psi_{\theta}$ (see \cite[Theorem 1.4.1]{De})

The general formulation of the Tamagawa number conjecture at the
non-critical values following \cite{Ka1} assumes $w+l+1>w$
\cite[Conjecture 2.2.7]{Ka1} because then one avoids the poles in
the bad Euler factors, and therefore the assumption $l\geq 0$.
But, for $M_{\theta}(w+l+1)$, there are no poles in the bad Euler
factors, see Remark \ref{rem26}. Thus, we can study the Tamagawa
number conjecture for $l<0$ using only the regulators maps.

In this section we construct elements in $K$-theory for
$M_{\theta}(w+l+1)$ with $0<-l\leq min(a_{\theta},b_{\theta})$ and
we study the image of these elements by the Beilinson regulator
map and the Soul\'e regulator map, obtaining the weak $p$-part of
the Tamagawa number conjecture.

Deninger \cite[pp.142-144]{De} already constructed elements in
$K$-theory for the motive $M_{\theta}(w+l+1)$ with $l<0$
non-critical and obtains their image by the Beilinson regulator
map, proving the Beilinson conjecture. He constructed these
elements in $K$-theory by use of a projector map
$\underline{\Kk}_{\Mm}$ without using complex multiplication. The
problem of his construction is that the Weil pairing appearing in
\S 5 to a $E[p^r]$-torsion point $\widetilde{t_r}$,
$\gamma(\widetilde{t_r})=<\widetilde{t_r},\widetilde{t_r}>$ is
trivial and the arguments through \S 5 does not generalize in
order to construct an Euler system to control the image by the
Soul\'e regulator map. We modify Deninger's projector map by
$\Kk_{\Mm}'$ (we use now complex multiplication), and we construct
the elements in $K$-theory using $\Kk_{\Mm}'$ and we reobtain
Beilinson's conjecture. With this modification the arguments in
the $p$-part of the weak Tamagawa number conjecture, i.e. the
image by the Soul\'e regulator map of these $K$-theory elements \S
4,\S 5, apply straightforward obtaining the weak $p$-part of the
Tamagawa number conjecture for $l<0$, Theorems \ref{lnegwTNCk},
\ref{lnegwTNCq}.

\end{subsection}

\begin{subsection}{Modification of Deninger's projector map.
Beilinson conjecture revisited.} Let us fix $w\geq 1$ and $l<0$
such that $-w-2l\leq-3$ with $0<-l\leq min(a_{\theta},b_{\theta})$
and let us consider the motive $M_{\theta}(w+l+1)$. With the fixed
embedding we have
$\vartheta=(\lambda_1,\ldots,\lambda_w)\in\theta_K$ and set
$I_1=\{i|\lambda_i\in Hom_K(K,\C)\}$ and $I_2=\{i|\lambda_i\notin
Hom_K(K,\C)\}$ and we have now that $0<|l|\leq \#I_1=a_{\theta}$
and $0<|l|\leq \#I_2=b_{\theta}$, where $|l|$ is the absolute
value. Denote by $\Delta=id^1\times id^2:E\rightarrow E\times E$
the diagonal map and by
$\Delta_{CM}=id^{1,CM}\times{id}^{2,CM}:E\rightarrow E\times E$
given by $e\mapsto(e,(\sqrt{d_K})e)$ where we understand
$\sqrt{d_K}\in End(E)$. Let us choose exactly $\#|l|$ elements in
the sets $I_1$ and $I_2$, denote their in increasing order
$i_1,\ldots,i_{|l|}\in I_1$ and $j_1,\ldots,j_{|l|}\in I_2$. Let
us define the projector map $pr:E^{w+l}\rightarrow E^{w+2l}$ by
the projection of the first $w+2l$-components of $E^{w+l}$ and
define $(id\times\Delta^{|l|}):E^{w+l}\rightarrow E^{w}$ (which it
depends of the choice in the sets $I_1$ and $I_2$) by
$(e_1,\ldots,e_{w+2l},e_{w+2l+1},\ldots,e_{w+l})\mapsto(e_{\alpha_1},\ldots,e_{\alpha_w})$
where $e_{\alpha_s}$ is defined as follows:

\begin{itemize}
\item if $\alpha_s$ appears in one component of the set of tuples
$L:=\{(i_1,j_1),\ldots, (i_{|l|},j_{|l|})\}$ then
$$e_{\alpha_s}:=\left\{\begin{array}{c}
id^{1}(e_{w+2l+m})\ if\ \alpha_s=i_m\\
id^2(e_{w+2l+m})\ if\ \alpha_s=j_m\\
\end{array}\right. ,$$
\item in the other case, then it is defined by
$e_{\alpha_s}:=e_{\tilde{n}}$ with $1\leq \tilde{n}\leq w+2l$ such
that $\alpha_s=\tilde{n}+\sum 1$ where the sum runs the naturals
that appear in some component of the elements of $L$ and which are
lower than $\alpha_s$.
\end{itemize}

We define the map $(id\times \Delta_{CM}^{|l|})$ similar as
$(id\times\Delta^{|l|})$ but replacing $id^i$ by $id^{i,CM}$.

The projector map $\Kk_{\Mm}'$ is defined by the commutative
diagram

\begin{small}
\begin{center}
\begin{tabular}{lcr}
$H^{2l+w+1}_{\Mm}(Sym^{2l+w}h^1E,\Q(w+2l+1))$&$\overset{pr^{*}}{\longrightarrow}$&$H^{2l+w+1}_{\Mm}(E^{2l+w+|l|},\Q(2l+w+1))$\\
$\mbox{}\ \mbox{} \ \mbox{} \Kk_{\Mm}'\downarrow$&&$\downarrow (id\times \Delta_{CM}^{|l|})_*\mbox{}\ \mbox{} \ \mbox{} $\\
$H^{w+1}_{\Mm}(M_{\theta\Q},\Q(w+l+1))$&$\overset{e_{\theta}}{\longleftarrow}$&$H^{w+1}_{\Mm}(h^1(E)^{\otimes w},\Q(l+w+1)).$\\
\end{tabular}
\end{center}
\end{small}

Deninger defines a projector map $\underline{\Kk}_{\Mm}$ with a
similar diagram as for our $\Kk_{\Mm}'$ but replacing the map
$(id\times\Delta_{CM}^{|l|})_{*}$ by the map
$(id\times\Delta^{|l|})_*$.

Let us choose the element in
$H^{w+1}_{\Mm}(M_{\theta\Q},\Q(w+l+1))$
$$\Upsilon_{\theta}:=\Kk_{\Mm}'\Ee_{\Mm}^{2l+w}(N_{K(E[\f])/K}((\Omega f^{-1}))),$$
where $\Ee_{\Mm}^{2l+w}$ is the Eisenstein symbol, $f$ a generator
of $\f_{\theta}$, $\Omega$ the period of $E$ and $(\Omega f^{-1})$
means the divisor in $\Z[E[\f_{\theta}]\setminus 0]$.

The next result is a modification of Deninger's result
\cite[pp.143-145]{De}.
\begin{thm} \label{beiapa} Suppose $a_{\theta}\not\equiv b_{\theta}$ mod $(\#\Y_K^*)$
with $a_{\theta}, b_{\theta}\geq 0$, $l<0$,
$w=a_{\theta}+b_{\theta}$, with $-w-2l\leq -3$ and $-l\leq
min(a_{\theta},b_{\theta})$. Define, up to sign,
$$\xi_{\theta,l}:=\frac{(\sqrt{d_K})^{2l}(2l+w)!L_p(\overline{\psi_{\theta}},-l)^{-1}
\Phi(\f_{\theta})}{2^{-1}N_{K/\Q}\f_{\theta}^{l}\psi_{\theta}(\rho_{\theta})\Phi(\f)}\Upsilon_{\theta}$$
which belongs to $H^{w+1}_{\Mm}(M_{\theta\Q},\Q(w+l+1))$ where
$L_p(\overline{\psi_{\theta}},-l)$ means the product of the Euler
factors of the primes above $p$ of $K$ at $-l$ (is well defined by
Remark \ref{rem26}), and $\rho_{\theta}$ is the id\`ele of $K$
such that $v_{\mathfrak{q}}(\rho_{\mathfrak{q}}^{-1}-f^{-1})\geq
0$ for $\mathfrak{q}\mid \f_{\theta}$ and
$v_{\mathfrak{q}}(\rho_{\mathfrak{q}})=0$ in the other primes
$\mathfrak{q}$. Then
$$r_{\Dd}(\xi_{\theta,l})=L_S^*(\overline{\psi_{\theta}},-l)\eta_{\theta},$$
where $S$ are the set of primes of $K$ that divide $\f_{\theta}p$,
$\eta_{\theta}$ is an $\Y_K[1/D_K]$-basis for
$H_B^w(M_{\theta\C},\Z(w+l))\otimes_{\Y_K}\Y_K[1/D_K]$ and
$L_S^*(\overline{\psi}_{\theta},-l)=\lim_{s+l\rightarrow
0}\frac{L_S(\overline{\psi}_{\theta},s)}{s+l}$.
\end{thm}
\begin{proof} We will follow closely Deninger's papers
\cite{Ded} and \cite{De}, we follow also in this proof his
notation where his $n$ is our $w+2l$. Deninger defines the element
$\xi_{\theta,l}$ from
$\underline{\Kk}_{\Mm}\Ee_{\Mm}^{2l+w}(N_{K(E[\f])/K}((\Omega
f^{-1})))$ instead of $\Upsilon_{\theta}$. We modify only the
calculation in \cite[(2.13)Lemma]{De} for $\Kk_{\Mm}'$ instead of
$\underline{\Kk}_{\Mm}$. One obtains (up to sign)
$$\frac{1}{(2\pi i)^w}\int_{E^w}\tilde{\Kk'_{\mathcal{D}}}(\tilde{\xi})\wedge dz^{(\underline{\varepsilon})}=$$
$$B_{\underline{\varepsilon}}\sqrt{d_K}^{|l|} \left(\begin{array}{c}
n\\
n+|l|-|\underline{\varepsilon}|\\
\end{array}\right)^{-1}
A(\Gamma)^{n+|l|}c_{n+|l|-|\underline{\varepsilon}|}$$ see the
calculation at the top of \cite[p.63]{Ded}. To precise the sgn we
should control the chosen order of the factors of the map
$(id\times\Delta_{CM}^{|l|})$, but for our interest is
unnecessary. Then the argument \cite[p.143-144]{De} applies in our
situation obtaining,
$$r_{\Dd}(\Upsilon_{\theta})=t_{\theta}L^*(\overline{\psi}_{\theta},-l)\eta_{\theta}$$
where $t_{\theta}$ is given by
$\frac{2^{-1}N_{K/\Q}\f_{\theta}^{l}\psi_{\theta}(\rho_{\theta})\Phi(\f)}{(\sqrt{d_K})^{2l}(2l+w)!
\Phi(\f_{\theta})}$(up to sign). By Remark \ref{rem26} we can
introduce the Euler factors above $p$ in the constant factor
$t_{\theta}$, obtaining the statement.
\end{proof}

\end{subsection}
\begin{subsection}{The weak Tamagawa number conjecture for $l<0$}

Following \S3 we define for $l<0$ the constructible module by
$$\Rr_{\theta}:=\xi_{\theta,l}\Y_K,$$
where $\xi_{\theta,l}$ is defined in theorem \ref{beiapa}. Let us
observe that with this notation we can follow straightforward all
the results and proofs of \S3 and \S4. In \S5 we need to compute
$\Kk_{\Mm}'\circ\Ee_{\Mm}^{w+2l}(N_{K(E[\f])/K}(\Omega f^{-1}))$.
We remember that we suppose once and for all that $p\nmid D_K$.
Denote by $e=(\widetilde{t_r})_r$ and element of the Tate module
$T_pE$ where $\widetilde{t_r}\in E[p^r]$ a $p^r$-torsion point for
$E$.
\begin{lem}\label{lema4.1} The
realization on Galois cohomology of the projector map $\Kk_{\Mm}'$
has the property, $\mathcal{K}_{\Mm}'(\widetilde{t_r}^{\otimes
2l+w})=e_{\theta}(\otimes^w\widetilde{t_r})\otimes
\gamma(\widetilde{t_r})^l$ where
$\gamma(\widetilde{t_r})=<\widetilde{t_r},\sqrt{d_K}\widetilde{t_r}>$.
\end{lem}
\begin{proof} Observe first that the projector map $\Kk_{\Mm}'$ is
$e_{\theta}\circ (id\times\Delta^{|l|}_{CM})_*\circ pr^*$. Let us
take $\delta_*:=(id\times\Delta^{|l|}_{CM})_*\circ pr^*$ and
observe that its transpose $\delta^*=pr_{*}\circ
(id\times\Delta^{|l|}_{CM})^*$ is part of the definition of
$\Kk_{\Mm}=e_{\theta}\circ \delta^*$ with $\underline{l}:=|l|>0$
given at \cite{K}. We want only to study these projector maps on
the Galois cohomology. Denote by $\mathcal{H}_{\Q_p}$ the \'etale
realization of $h^1(E)(1)$ and observe that there is an
isomorphism $\mathcal{H}_{\Q_p}^{*}(1)\cong \mathcal{H}_{\Q_p}$,
since $(h^1(E)(1))^*=h_1(E)(-1)\cong h^1(E)(1)(-1)=h^1(E)$. The
map $\delta^*$ is given by
$$H^1(\Y_S,Sym^{2\underline{l}+\underline{w}}(\mathcal{H}_{\Q_p}) (1))\rightarrow H^1(\Y_S,
Sym^{\underline{w}}(\mathcal{H}_{\Q_p})(\underline{l}+1)),$$ and
because the map $\delta_*$ is the transpose for the map
$\delta^*$, up to Tate twist by $\underline{w}+\underline{l}$, it
is represented by global Tate duality by,
$$H^1(\Y_S,
Sym^{\underline{w}}(\mathcal{H}_{\Q_p})^*(-\underline{l}-1)(1))\rightarrow
H^1(\Y_S,Sym^{2\underline{l}+\underline{w}}(\mathcal{H}_{\Q_p})^*(-1)
(1)).$$

Is known \cite{K} that

$$\delta^*(\plim{r}(\otimes^{2\underline{l}+\underline{w}} \widetilde{t}_r))=\plim{r}((\otimes^{\underline{w}} \widetilde{t}_r)
\gamma(\widetilde{t}_r)^{\underline{l}})$$ write this equality
also by $\delta^*(\otimes^{\underline{w}+2\underline{l}}
v)=(\otimes^{\underline{w}}v)\gamma(v)^{\underline{l}}$. Take now
the dual map by $Hom(\mbox{},\Z_p)$ and with the identification
$T_pE\cong Hom(T_pE,\Z_p(1))$, we obtain

$$(\otimes^{\underline{w}}v(-1))\gamma(v)^{-\underline{l}}\mapsto
\otimes^{\underline{w}+2\underline{l}}v(-1)$$ twisting now by
$\underline{w}+\underline{l}$ we arrive to the definition for
$\delta_*$ and,
$$\delta_*(\otimes^{\underline{w}} v)\mapsto (\otimes^{\underline{w}+2\underline{l}}v)\gamma(v)^{-\underline{l}}.$$
Now take this equality at level $r$,
$w=\underline{w}+2\underline{l}$, $l=-\underline{l}$, and apply
the idempotent $e_{\theta}$ to finish.
\end{proof}

 After
the lemma \ref{lema4.1} all the results of \S5 and the proofs of
\S5
 follow straightforward up to a power of 2 and $D_K$, (the reader could
make these modifications which follow only from our definition of
$\Rr_{\theta}$). Therefore we obtain the weak $p$-part of the
Tamagawa number conjecture with $K$-coefficients and
$\Q$-coefficients, under standard hypothesis from Iwasawa theory
for imaginary quadratic fields:

{\bf (***)} Let $p$ be a fix prime such that $p\nmid 6N_{K/\Q}\f$.
Suppose that $\psi_{\theta}$ has infinity type
$(a_{\theta},b_{\theta})$ with $a_{\theta}, b_{\theta}$
non-negative integers, such that $a_{\theta}\not\equiv b_{\theta}\
mod(\# \Y_K^*)$ and $w=a_{\theta}+b_{\theta}\ge1$ verifies
$-w-2l\le-3$ with $l<0$ and $-l\leq min(a_{\theta},b_{\theta})$.
Suppose that $\Y_K^*\rightarrow(\Y_K/\f_{\theta})^*$ is injective.
Suppose moreover that the representation $\chi$ of
$Gal(K(E[p])/K)$ in
$Hom_{\Y_p}(H^w(M_{\theta}\times_K\overline{K},\Z_p(w+l)),\Y_p)$
is a good representation which is not equal as
$\Delta$-representation to the cyclotomic representation.

\begin{thm}\label{lnegwTNCk} Assume hypotheses {\rm (***)}.
Then, there is an $\Y_K$-submodule $\Rr_{\theta}\subset H_{\Mm}$
of rank 1 such that:
\begin{enumerate}
\item
$det_{\Y_{K}[1/D_K]}(r_{\Dd}(\Rr_{\theta}\otimes_{\Y_K}\Y_K[1/D_K]))\cong$
$$L_S^*(\overline{\psi}_{\theta},-l)
det_{\Y_K[1/D_K]}(H^w_B(M_{\theta\C},\Z(w+l))\otimes_{\Y_K}\Y_K[1/D_K])$$
in
$det_{\Y_K[1/D_K]\otimes\R}(H_B^w(M_{\theta\C},\Z(w+l))\otimes_{\Y_K}\Y_K[1/D_K]\otimes\R)$.
\item The map $r_p$ induces an isomorphism
$$det_{\Y_{K}\otimes\Z_p}(\Rr_{\theta})\cong det_{\Y_{K}\otimes\Z_p}(R\Gamma(\Y_{K}[1/S],M_{\theta\Z_p}(w+l+1))^{-1}.$$
\end{enumerate}
Here
$$L_S^*(\overline{\psi}_{\theta},-l)=\lim_{s\rightarrow-l}\frac{L_S(\overline{\psi}_{\theta},s)}{s+l},$$
and $S$ is the set of primes of $K$ dividing $p$ and the ones
dividing $\f_{\theta}$.

Moreover, if $r_p$ is injective on $\Rr_{\theta}$, the second part
can be written as
$$det_{\Y_{K}\otimes\Z_p}(H^1(\Y_{K}[1/S],M_{\theta\Z_p}(w+l+1))/r_p(\Rr_{\theta}))\cong $$
$$det_{\Y_{K}\otimes\Z_p} H^2(\Y_{K}[1/S],M_{\theta\Z_p}(w+l+1)).$$

\end{thm}
\begin{thm}\label{lnegwTNCq}
Suppose hypotheses {\rm (***)}.

Then, there is a $\Z$-submodule $\Rr_{\theta}$ in $H_{\Mm}$ of
rank 2 such that:
\begin{enumerate}
\item The map $r_{\Dd}\otimes\R$ is an isomorphism restricted to
$\Rr_{\theta}\otimes\R$. \item
$dim_{\Q}(H^w_B(M_{\theta\C},\Q(w+l)))=ord_{s=-l}L_S(M_{\theta\Q},s)=2$.
\item We have the equality
$$r_{\Dd}(det_{\Z[1/D_K]}(\Rr_{\theta}\otimes_{\Y_K}\Y_K[1/D_K]))=$$
$$L_S^*(M_{\theta\Q},-l)det_{\Z[1/D_K]}(H^w_B(M_{\theta},\Z(w+l))\otimes_{\Y_K}\Y_K[1/D_K])$$
where
$$L^*_S(M_{\theta\Q},-l)=\lim_{s\rightarrow-l}\frac{L_S(M_{\theta\Q},s)}{(s+l)^2}$$
and $S$ is the set of places of $K$ that divides $p$ and the
places dividing the conductor $\f_{\theta}$. \item We have that
$$det_{\Z_p}(\Rr_{\theta}\otimes \Z_p)=det_{\Z_p}(R\Gamma(\Y_{K}[1/S],M_{\theta\Z_p}(w+l+1)))^{-1}.$$
If $r_p$ is injective on $\Rr_{\theta}$, then
$r_p(det_{\Z_p}(\Rr_{\theta}\otimes\Z_p))$ is a basis of the
$\Z_p$-lattice
$$det_{\Z_p}(R\Gamma(\Y_{K}[1/S],M_{\theta\Z_p}(w+l+1)))^{-1}$$
$$\subset det_{\Q_p}(R\Gamma(\Y_{K}[1/S],M_{\theta\Z_p}(w+l+1)\otimes\Q)[-1]).$$
\end{enumerate}
\end{thm}
\end{subsection}

\begin{section}{Some explicit examples}

Observe first if we consider the Hecke character associated to the
idempotent $e_{\theta}$ with infinite type $(1,0)$, then, Theorem
\ref{kht:teo} is exactly \cite[Theorem 1.1.5]{K} (hypothesis
($\diamond\diamond\diamond$) for the infinite type $(1,0)$
coincides with the ones that appears in \cite[Theorem 1.1.5]{K} ,
use \cite[Corollary 2.2.11]{K}, Lemmata \ref{4.5}, \ref{4.6},
\ref{5.5} and comments after Definition \ref{4.3}).

Let us give examples of Hecke characters of infinite type
different to $(1,0)$. In the following, take $E$ the elliptic
curve $y^2=4x^3-4x$ defined over $K=\Q(i)$ and the $e_{\theta}$'s
are defined from this fixed $E$. $E$ is an elliptic curve with CM
by $\Z[i]$, and for $p\geq 5$ we have $p\nmid
6N_{K/\Q}\mathfrak{f}$.

Let us consider $e_{\theta}$ with infinite type $(a,b)$ satisfying
$a\not\equiv b(mod\ 2)$, $a>b\geq 0$ and $-(a+b)-2l\leq -3$ with
$-l\leq b$. Take $p$ such that it splits in $\Q(i)$ and
$p-1>max(3,a-b)$. We have that $e_{\theta}$ and $p$ satisfy all
the hypothesis ($\diamond\diamond\diamond$) and $(***)$, with the
exception of the condition (A) in Definition \ref{4.3}, (use Lemma
\ref{4.6} and with a similar proof done for Lemma \ref{5.5} one
obtains that $(\Z[i])^*\rightarrow
(\Z[i]/\mathfrak{f}_{\theta})^*$ is injective if
$(\#\mathcal{O}_K^*,a-b)=1$). We impose $(a-b,p-1)=1$ to ensure
that $e_{\theta}$ and $p$ satisfy the condition (A) (see Lemma
\ref{4.5}).
For such $e_{\theta}$ and $p$ we obtain the conclusion of theorems
\ref{kht:teo}, \ref{qht:teo}, \ref{lnegwTNCk} and \ref{lnegwTNCq}.

For an explicit example take $e_{\theta_1}$ with infinite type
$(3,0)$ and $p=5$, then all hypothesis
($\diamond\diamond\diamond$) are satisfied and moreover $p=5$ is a
regular prime for $\Q(i)$ (see \cite[p.33]{Ya}) therefore by
theorem \ref{qht:teo} (and remark \ref{5.14}) we get the weak
5-part of the Tamagawa number conjecture for the dual of the
motive $M_{\theta_1}(3+l+1)$ twisted by 1 (the special value for
the motive $M$ is the special value $L(M,0)$, and in our
formulation we get the special value
$L(\overline{\psi}_{\theta},-l)=L((\tiny{M_{\theta}(w+l+1)})\check{\mbox{}}\
 (1),0)$). By use of the functional equation between the motive
and its dual twisted by 1 and good compatibilities, one should
obtain the 5-part of the Tamagawa number conjecture for the
motive $M_{\theta_1}(3+l+1)$ for any $l\geq 0$.\\
 For another
explicit example, take $e_{\theta_2}$ with infinity type
$(a,b)=(3,2)$ and $p=5$, and we get the weak 5-part of the
Tamagawa number conjecture for
$({M_{\theta_2}(3+2+l+1)})\check{\mbox{}}\ (1)$ for $l\geq -1$ by
theorems \ref{qht:teo} and \ref{lnegwTNCq}.

\end{section}

\begin{section}*{Acknowledgement}
I am indebt with X.Xarles. I would like to thank him for his
useful comments and suggestions and moreover to give me energy to
try to make me enjoy another time doing mathematical research. It
is also a big pleasure to thank G. Kings for many discussions, C.
Deninger for introducing me in the $L$-function world and K. Rubin
for clarifying a doubt on elliptic units. Finally I thank the
referees for their comments and suggestions.
\end{section}

Francesc Bars Cortina\\
Depart. Matem\`atiques, Edifici C,\\
 Universitat Aut\`onoma de
Barcelona. \\
08193 Bellaterra.\\
 Catalonia. Spain\\
 E-mail:~\textsf{francesc@mat.uab.cat},

\end{document}